\documentclass[runningheads]{llncs}

\usepackage{amsmath,amsfonts}
\usepackage{amssymb}

\usepackage[final]{graphicx}
\usepackage{url}
\newcommand{\keywords}[1]{\par\addvspace\baselineskip\noindent\keywordname\enspace\ignorespaces#1}

\usepackage{color}
\usepackage{algorithm2e}
\usepackage{algorithmic}
\usepackage{subfigure}

\newlength\myindent
\begin{document}


\title{
Stochastic Image Deformation in Frequency Domain and Parameter Estimation using Moment Evolutions
}

\author{Line K\"uhnel\inst{1}, Alexis Arnaudon\inst{2}, Tom Fletcher\inst{3} and Stefan Sommer\inst{1}}
\institute{
Department of Computer Science (DIKU), University of Copenhagen, Denmark
\and 
Department of Mathematics Imperial College London, UK
\and
Department of Electrical and Computer Engineering / Department of Computer Science,
University of Virginia, USA
}

\maketitle

\begin{abstract}
Modelling deformation of anatomical objects observed in medical images can help describe disease progression patterns and variations in anatomy across populations. We apply a stochastic generalisation of the Large Deformation Diffeomorphic Metric Mapping (LDDMM) framework to model differences in the evolution of anatomical objects detected in populations of image data. The computational challenges that are prevalent even in the deterministic LDDMM setting are handled by extending the FLASH LDDMM representation to the stochastic setting keeping a finite discretisation of the infinite dimensional space of image deformations. In this computationally efficient setting, we perform estimation to infer parameters for noise correlations and local variability in datasets of images. Fundamental for the optimisation procedure is using the finite dimensional Fourier representation to derive approximations of the evolution of moments for the stochastic warps. Particularly, the first moment allows us to infer deformation mean trajectories. The second moment encodes variation around the mean, and thus provides information on the noise correlation. We show on simulated datasets of 2D MR brain images that the estimation algorithm can successfully recover parameters of the stochastic model.
\end{abstract}

\section{Introduction}
\label{sec:int}

Classical models describing the evolution of anatomical objects, occurring from child development, from natural ageing, or from disease processes, are generally smooth and deterministic. However,
when analysing such deformations across populations of subjects, the individual deviations have to be taken into account as they otherwise affect the average trend in the evolution of the population. It is natural to assume that the subject-specific deformations are not purely deterministic and that stochastic variation may occur at any time of the evolution process. In this paper, we develop the technical framework to model such combinations of population average and subject-specific stochastic evolutions. 
\begin{figure}[t]
    \begin{center}
        \includegraphics[scale=0.8]{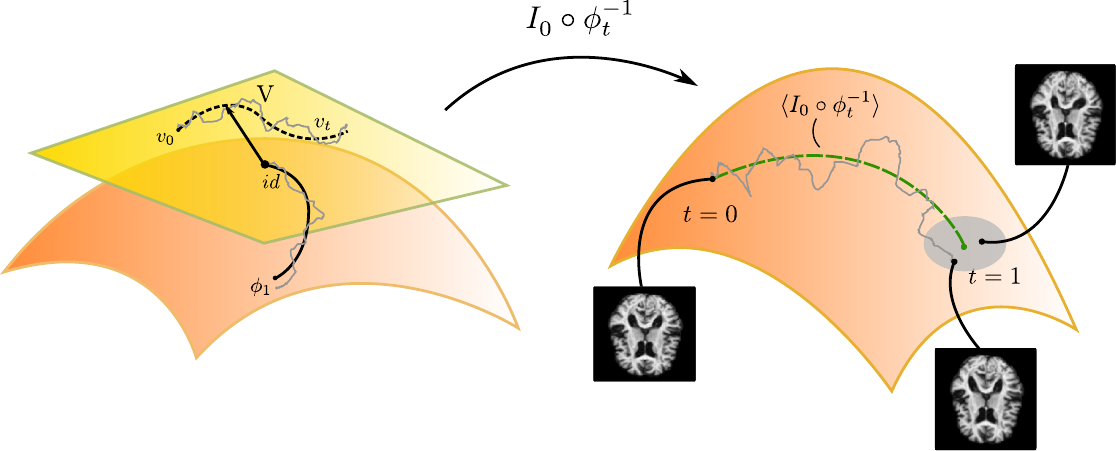}
    \end{center}
    \caption{The presented stochastic deformation model includes both the average population trend over time starting at time $t=0$, and the subject specific stochastic evolution.
      The deformation $\phi^{-1}_t$ is modelled as a stochastic perturbation of the LDDMM flow of diffeomorphisms resulting in a stochastic deformation $I_0\circ\phi^{-1}_t$ of the image $I_0$. 
      The stochastic deformation forms at each time point a distribution around the population trend describing uncertainty in the evolution. 
      The mean flow $\langle I_0\circ\phi^{-1}_t\rangle$ (dashed green line) models the general population evolution.}
    \label{fig:model}
\end{figure}

Shape changes of anatomical objects, observed in an image $I_0$, can be modelled using image registration. Here, the goal is to determine a deformation $\phi\colon\mathcal{D}\to\mathcal{D}$ of the image domain $\mathcal{D}$, minimising an energy
    $E(\phi) = R(\phi) + S(I_0\circ\phi^{-1},I_1)\,,$
for a regularisation $R$ and similarity measure $S$, see for example \cite{younes2010shapes}. 
In this framework, the deformation $\phi$ is a deterministic, smooth, bijective function independent of time. 
The Large Deformation Diffeomorphic Metric Mapping framework (LDDMM, see e.g. \cite{beg_computing_2005}), defines the deformation $\phi$ as the endpoint of a flow of diffeomorphisms $\phi_t$ solution to the differential equation
\begin{align}
  \frac{d}{dt}\phi_t(x,y) = v_t\circ\phi_t(x,y)\,,
  \label{reconstruction-det}
\end{align}
with initial value $\phi_0 = \text{id}$.
The time-dependent velocity field $v_t$ is for each $t$ an element of the space of vector fields on $\mathcal{D}$, $V=\mathfrak{X}(D)$, and solves the EPDiff equation $\frac{d}{dt}v_t = -K\text{ad}^*_{v_t}m_t$. 
Here $m_t = \frac{\delta l}{\delta v_t} = K^{-1} v_t = Lv_t$ is the conjugate momentum of the velocity $v_t$ if the regularisation $R(\phi)= l(v) = \frac12 \|v\|^2_K$ corresponds to the kinetic energy of the flow with a reproducible kernel $K$, see for example \cite{younes2010shapes} and references therein.
Considering deformations via a flow in the diffeomorphism group introduces a natural time component which can be used to model the evolution of anatomical objects over time. 

The flow of diffeomorphisms in the LDDMM model is deterministic, and it is hence only possible to introduce uncertainties in the initial velocity field $v_0$. 
Random variation of the initial velocity field has been discussed in, for example, the random orbit model of~\cite{miller_statistical_1997} and, for Bayesian principal geodesic analysis, \cite{zhang_bayesian_2015}.

Modelling uncertainty by a random initial velocity field implies that the entire variation over time is the result of uncertainty at the initial point of the evolution process. For longitudinal models, it is arguably more natural to model uncertainty time-continuously and thus having randomness occurring through the entire evolution process.
To enable random evolution, we aim at modelling the variation as the endpoint of a \emph{stochastic} flow of deformations, $\phi_t$. This, in turn, gives a separation between the deterministic mean population evolution, and the per subject individual stochastic trajectories as illustrated in Fig.~\ref{fig:model}. 

In this work, we use the stochastic LDDMM framework presented in~\cite{arnaudon_geometric_2017}, based on the stochastic fluid dynamics model introduced in \cite{holm_variational_2015}. 
Arnaudon et al.~\cite{arnaudon_geometric_2017} used a stochastic flow of diffeomorphisms to model uncertainties in the evolution of any data type on which the diffeomorphism group acts, in particular for landmarks.
The stochastic flow is defined by a stochastic differential equation (SDE) with a diffusion term parametrised by Eulerian noise fields $\sigma(x)$ instead of a standard Lagrangian noise associated with the flow.
In the latter case, \cite{TrVi2012,vialard2013extension} studied a stochastic model of landmarks dynamics with a different noise for each landmark and more recently \cite{marsland2016langevin} added dissipation to this model. 

The stochastic LDDMM with Eulerian noise fields of \cite{arnaudon_geometric_2017} applies to any data structure on which the diffeomorphism group acts without modification of the noise structure. Thus, in addition to landmarks, which have a finite dimensional structure, the framework can be applied to complete images \cite{arnaudon_string_2018}, provided that a good finite dimensional approximation can be applied to make numerical simulations possible. 
In this paper, we extend the fast LDDMM solver FLASH presented in~\cite{zhang2015finite} to include the stochastic deformations.
The algorithm is based on spatial Fourier transformation of a stochastic version of the EPDiff equation which results in a natural truncation of the high-frequencies in the stochastic process $v_t$. 
For Eulerian noise field, this is not an issue as we have total control on the spatial correlation of the noise (and white noise in time). Thus we avoid introducing high frequencies with the noise and thus obtain a good numerical approximation of the flow.
This truncation gives rise to a dimensionality reduction resulting in a significant computational speed up, and it makes the use of noise in image matching possible and relevant for applications. 

With the paper, we make the following contributions: 1) we incorporate stochasticity in the FLASH framework to model stochastic image evolutions in a finite-dimensional setting; 2) we derive the Fokker-Planck equations in the finite dimensional model and use this to approximate the evolution of the moments of the stochastic flow and deformed images; 3) we show how matching of the moment images can be used to estimate unknown parameters in the model; 4) we illustrate the use of the model on 2D brain images by recovering parameters from images simulated from the model.

\section{Stochastic Image Deformation}
\label{sec:stocev}

The LDDMM framework models deformation over a time region as a flow, $\phi_t$, in the space of diffeomorphisms $\text{Diff}(\mathcal{D})$. As $\phi_t$ varies smoothly and deterministically, applying the flow to an image $I_0$ results in a time evolution of the image which does not describe any uncertainty in the deformation. In this section, we describe the stochastic LDDMM extension \cite{arnaudon_geometric_2017} which exactly models the uncertainty of deformations.

\begin{figure}[htb!]
    \begin{center}
        \subfigure{\includegraphics[width=0.24\columnwidth,trim=20 20 35 20,clip]{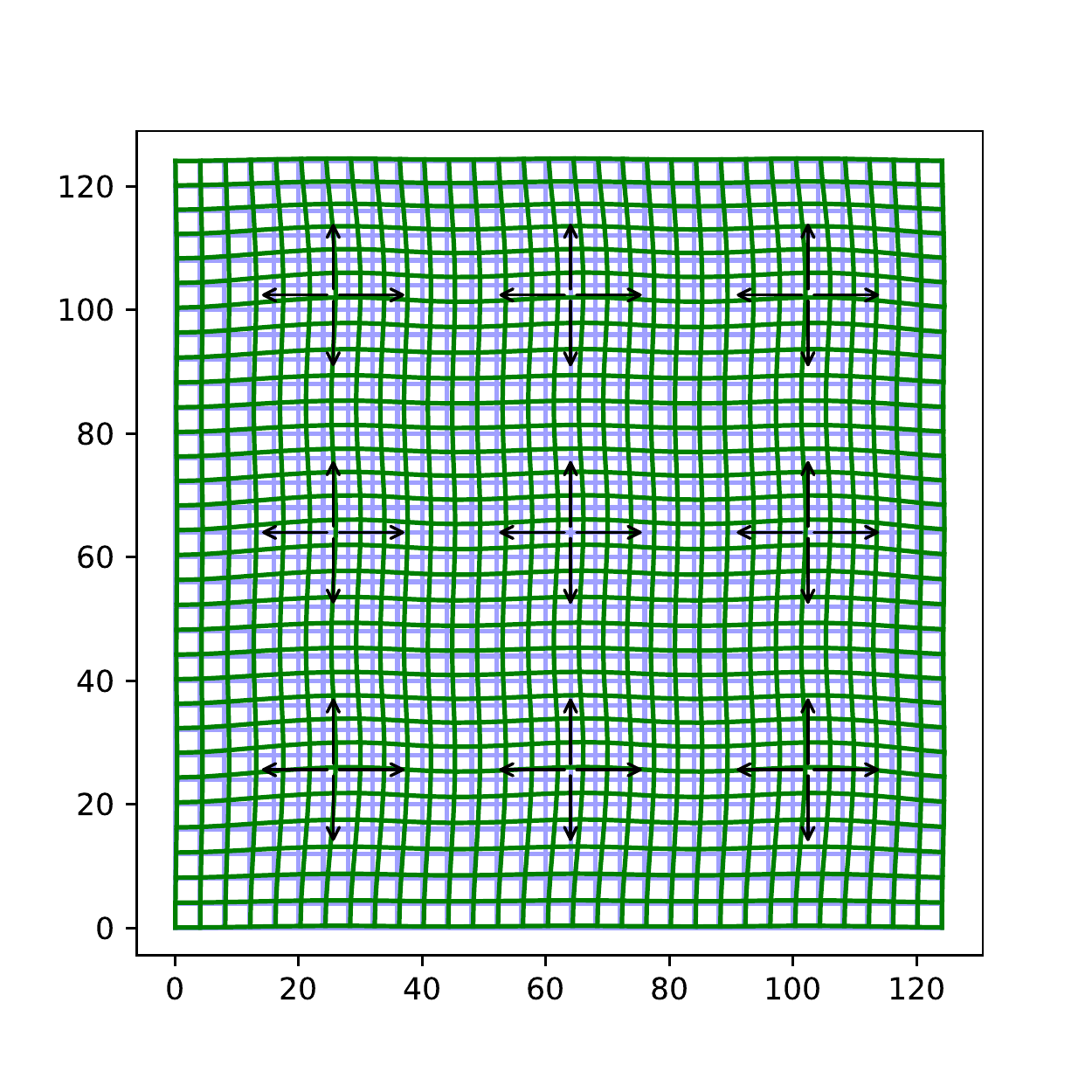}}
        \subfigure{\includegraphics[width=0.24\columnwidth,trim=100 32 100 30,clip,angle=90]{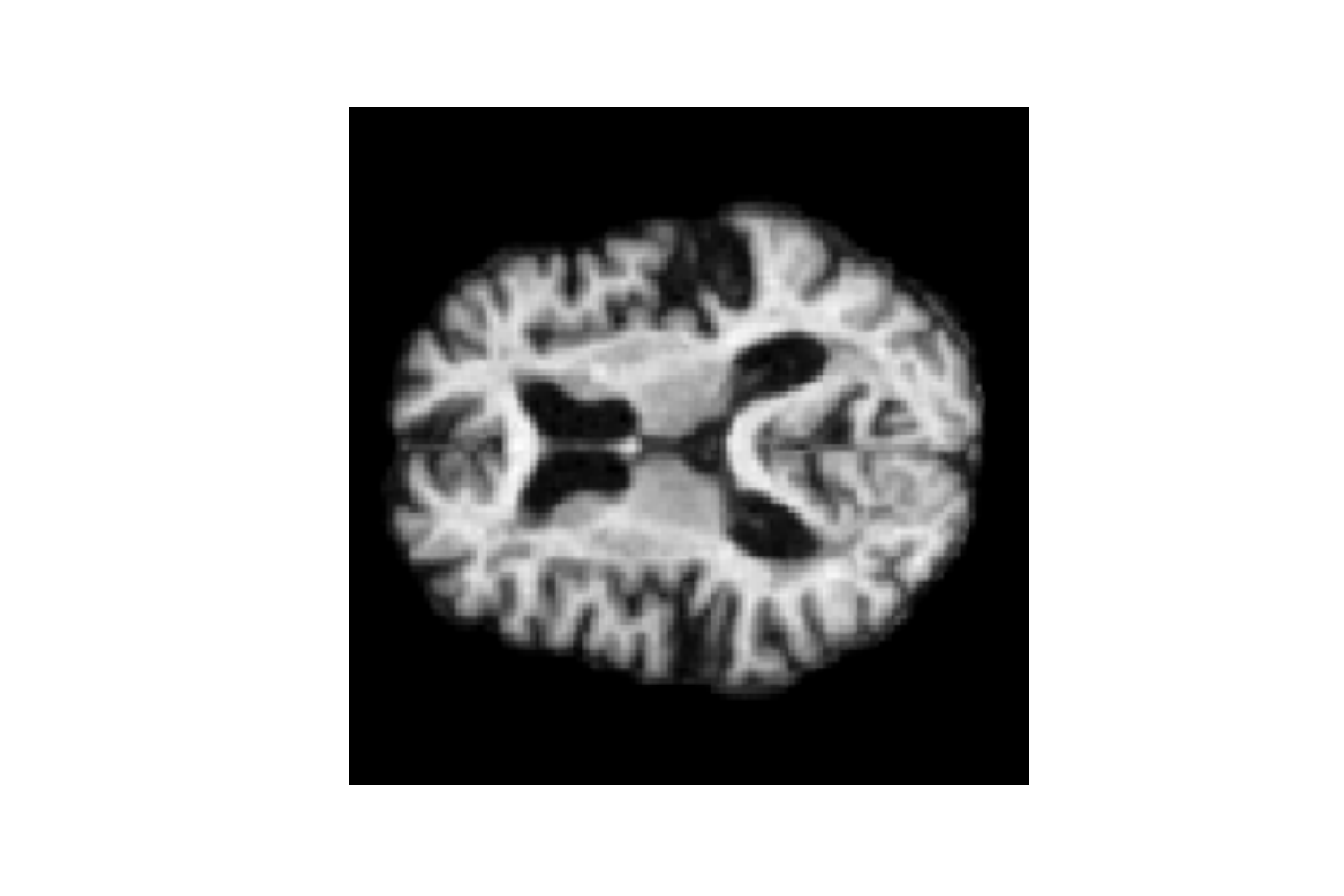}}
        \subfigure{\includegraphics[width=0.24\columnwidth,trim=100 30 100 30,clip,angle=90]{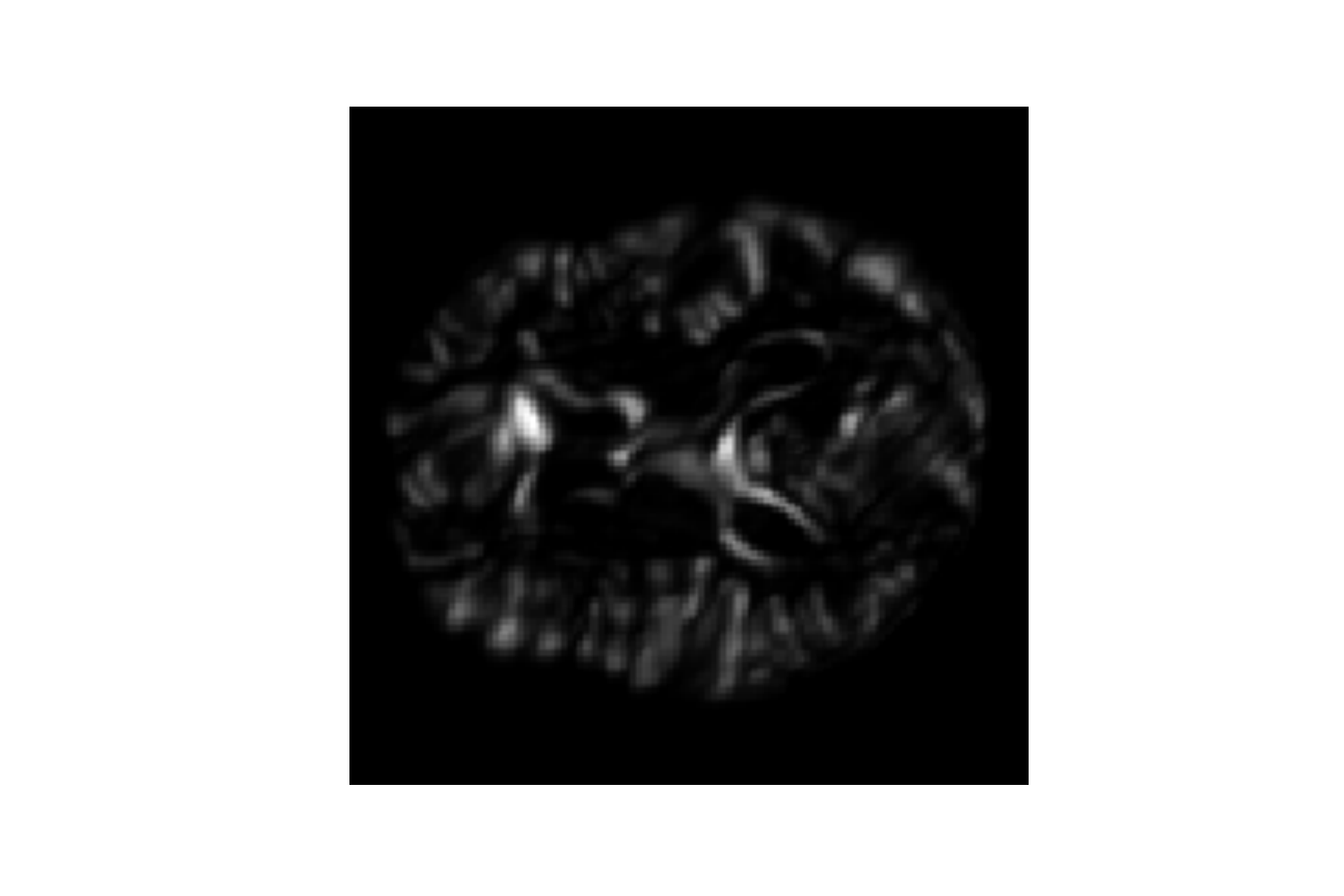}}\vspace{-0.2cm}\\
        \subfigure{\includegraphics[width=0.24\columnwidth,trim=100 30 100 30,clip,angle=90]{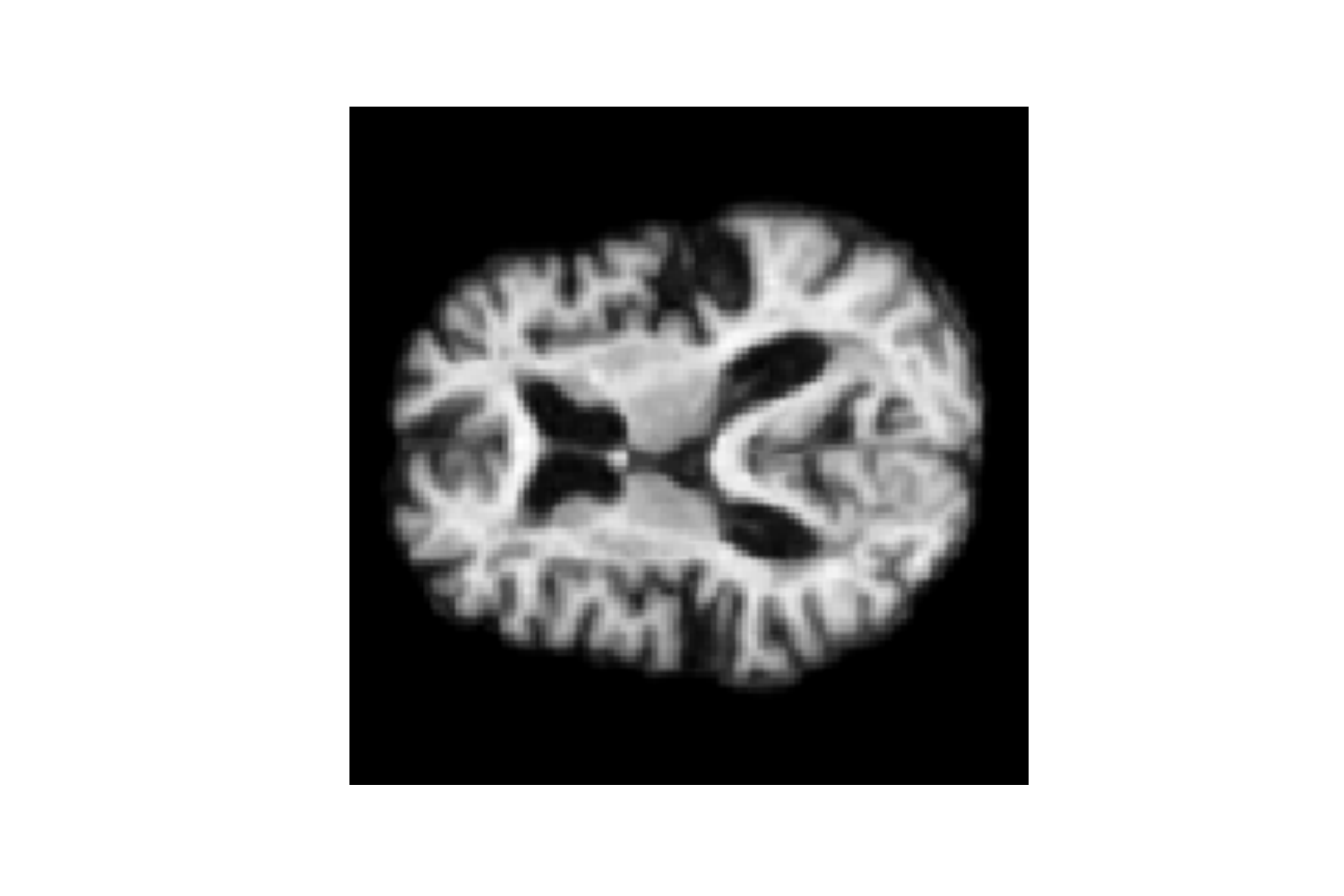}}
        \subfigure{\includegraphics[width=0.24\columnwidth,trim=100 30 100 30,clip,angle=90]{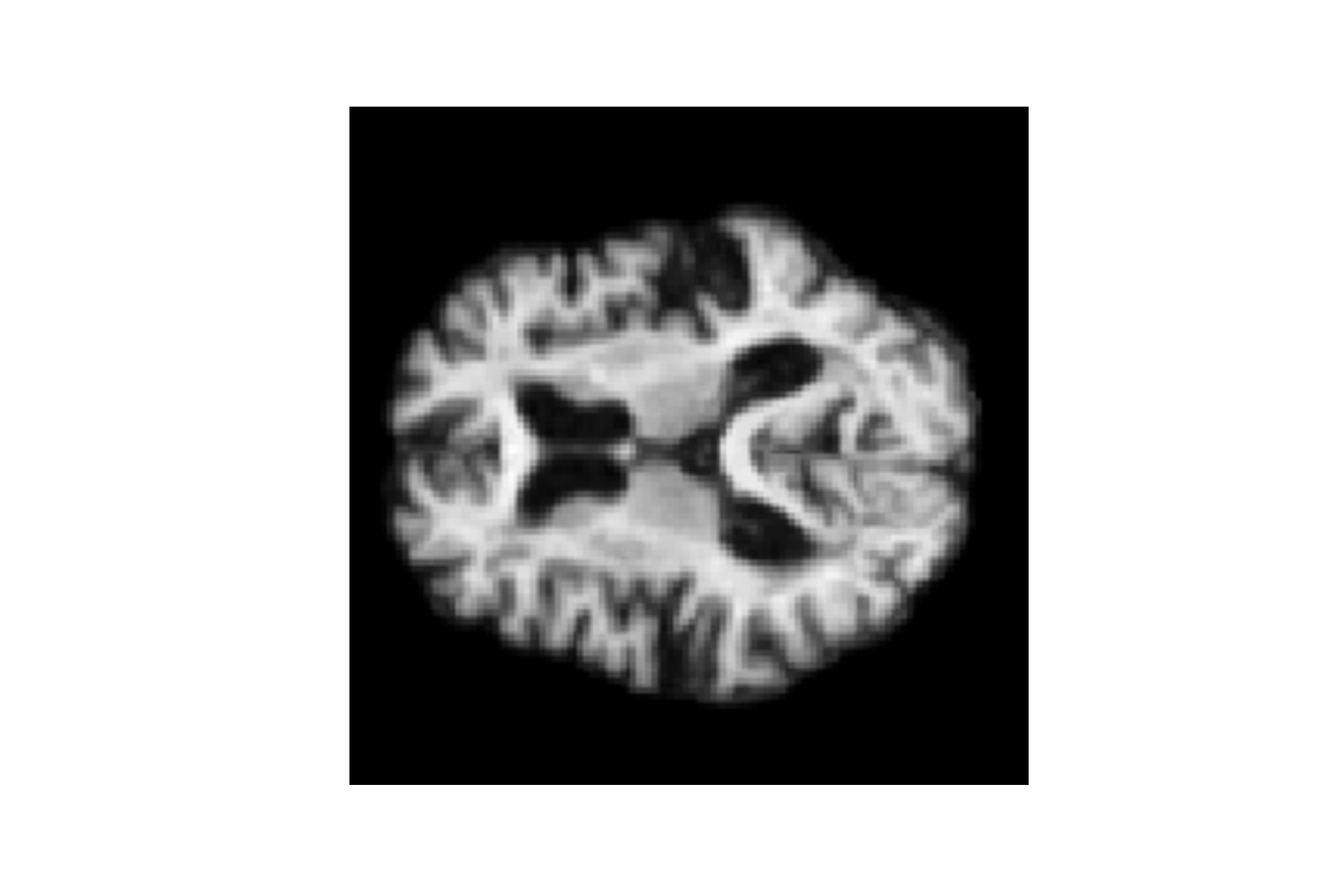}}
        \subfigure{\includegraphics[width=0.24\columnwidth,trim=100 30 100 30,clip,angle=90]{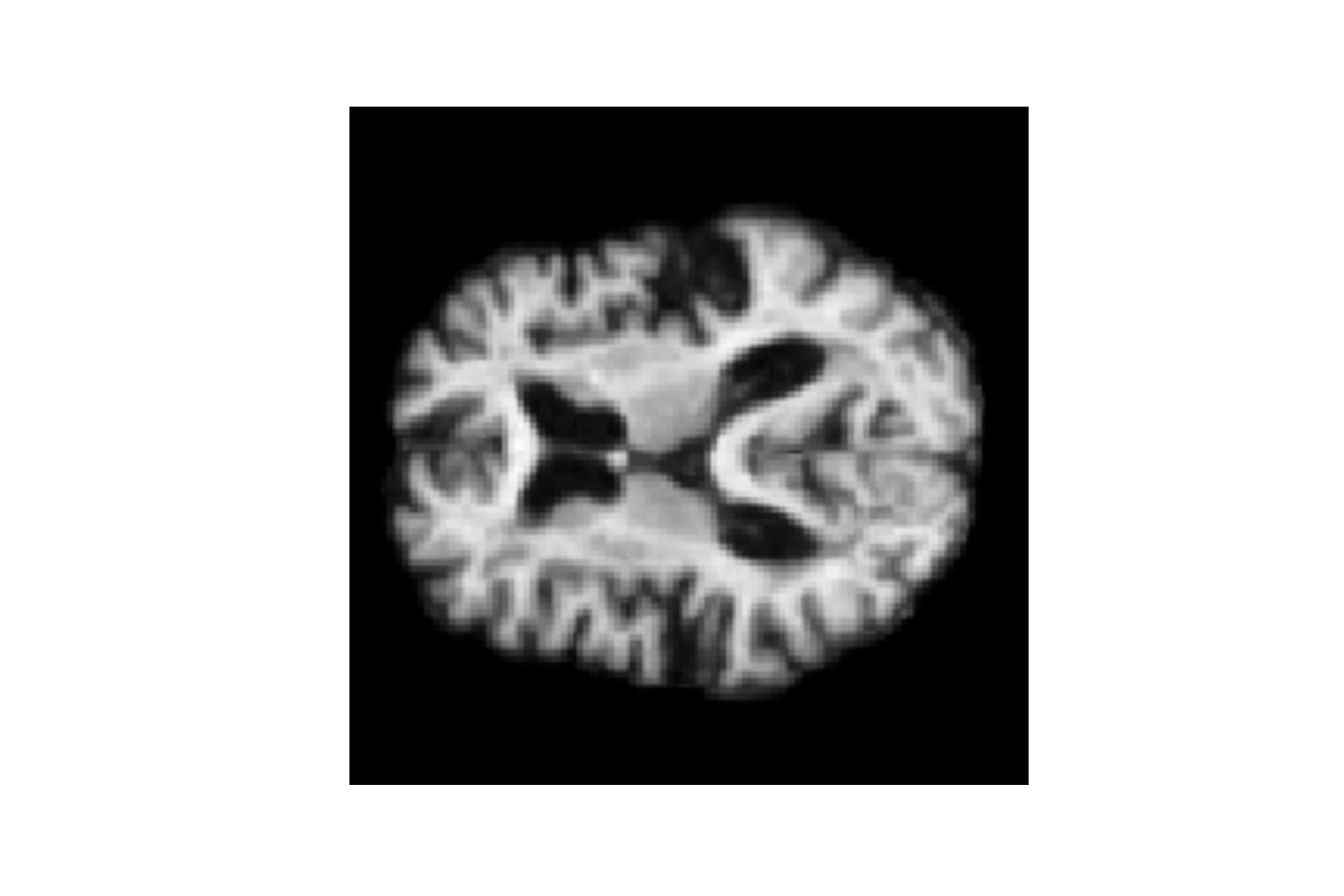}}
    \caption{\textit{First row:} (left) The noise fields location on the domain. (middle) Initial brain, $I_0$. (right) Variation in the data sample. \textit{Second row:} 3 examples of simulated data by the stochastic deformation $\phi_t^{-1}$.}
    \label{fig:init_var}
    \end{center}
\end{figure}

We consider a probability space $(\Omega,P,\mathcal{F})$ and let $\phi_t\colon\Omega\times [0,1]\times\mathcal{D}\to\mathcal{D}$ be a stochastic process on the time interval $[0,1]$, i.e. for each $t\in [0,1], \ \omega\in\Omega$, $\phi_t(\omega)\colon\mathcal{D}\to\mathcal{D}$ is a deterministic deformation of the image domain $\mathcal{D}$. 
The stochastic flow $\phi_t$ is defined as a Stratonovich SDE with diffusion term based on noise fields, $\sigma_k\in\mathfrak{X}(\mathcal{D})$, $k = 1,\ldots,p$, on the image domain $\mathcal{D}$. See Fig. \ref{fig:init_var} for an example of noise fields describing the stochastic flow of deformations. The Stratonovich SDE defining the flow $\phi_t$ is,
\begin{align}
  d\phi_t = v_t(\phi_t) dt + \sum_{k=1}^p\sigma_k(\phi_t)\circ_S dW_t^k, \quad \phi_0 = \mathrm{id}\, ,
\label{eq:def}
\end{align}
where $W^k_t$ denotes $p$ one-dimensional Wiener processes (or Brownian motions) on $\mathbb{R}$ and $\circ_S$ Stratonovich integration. 
The time-varying velocity field $v_t$ is the solution of the stochastic Euler-Poincar\'e equation,
\begin{align}
  dv_t = -K\mathrm{ad}_{v_t}^* m_t - \sum_{k=1}^p K\mathrm{ad}_{\sigma_k}^* m_t\circ_S dW^k_t, \quad v_0\in\mathfrak{X}(\mathcal{D})\, ,
\label{eq:vt}
\end{align}
for the same momentum $m_t = Lv_t$ as in the deterministic equations. 
In~\cite{arnaudon_geometric_2017}, it was shown that under the stochastic deformation \eqref{eq:def}, the momentum is preserved and \eqref{eq:vt} is hence the stochastic version of the EPDiff equation. 

Zhang et al.~\cite{zhang2015finite} observed that the last operation of the EPDiff equation is applying a low-pass filter $K$. 
As $K$ is an operator suppressing all high-frequencies in Fourier domain, performing dimensionality reduction on the number of frequencies results in a large computational gain, and, importantly, only a restricted amount of information being removed. Similarly to the drift, the stochastic term of \eqref{eq:vt} has the smoothing operator $K$, applied exclusively to the spatially dependent noise amplitude, and not $dW_t^k$. Hence, it is also possible to benefit from the computational speedup of FLASH in the stochastic setting. We make a spatial Fourier transformation of $dv_t$ and truncate the high-frequencies of this stochastic field. Details on the FLASH Fourier space calculation in the deterministic setting can be found in \cite{zhang2015finite}. As an example, stochastic shooting of an $128\times 128$ image truncated to $16$ Fourier frequencies with $100$ times steps and $1$ noise field on a standard laptop (i7 processor) takes approximately $1.5$ seconds.

The stochastic deformation $\phi_t$ models both the population trend and the subject-specific variations. The population trend is represented by the deterministic part of \eqref{eq:def} and \eqref{eq:vt}. It is a function of the initial velocity field $v_0$ and describes the global trend of the population, e.g. ageing or a disease progression. The noise fields, on the other hand, describe the subject-specific variation present in the data. 
The noise fields are modelled as local objects of the full dimensional data space with correlations decreasing with the distance to the centre of the object, e.g. Gaussian kernels on the domain.

\section{Moment Matching}\label{sec:mm}

Consider $n$ individuals observed at $t=0$, $I_0^1,\ldots,I_0^n$ and again at $t=1$, $I_1^1,\ldots,I_1^n$. The goal is to model the population trend and infer the noise structure describing the variation in the observed data at $t=1$. 
 For this, we aim at estimating the parameters of the stochastic deformation model. These parameters are the noise fields $\sigma_1,\ldots,\sigma_p$, the initial velocity field $v_0$, and parameters of the LDDMM RKHS kernel $K$. In this paper, we focus on estimating the parameters of the noise fields $\sigma_1,\ldots,\sigma_p$.

We base parameter estimation on method of moments, i.e. we seek to match moments of the observed data distribution with moments of the distribution of images at $t=1$ generated by the stochastic deformation model. 
For simplicity, we remove the subject effect at $t=0$ by considering a single initial image $I_0$, the average $\hat{I}_0=\frac{1}{n}\sum_{i=1}^n I_0^i$. 
This implies that population variation at time $t=0$ is removed and that the model, therefore, needs to account for the entire population variation at $t=1$. See Fig. \ref{fig:model} for a visualisation of the model.
Moments of the random variable $\hat{I}_0\circ\phi_1^{-1}$ are matched to the data moments, to retrieve the parameters for the noise fields $\sigma_k$.
 The model can handle subject specific initial images with minor changes to the moment equations presented in Section \ref{sec:taylor}.

\begin{figure}[htb!]
    \begin{center}
        \subfigure[$t=0$]{\includegraphics[width=0.24\columnwidth,trim=100 30 100 30,clip,angle=90]{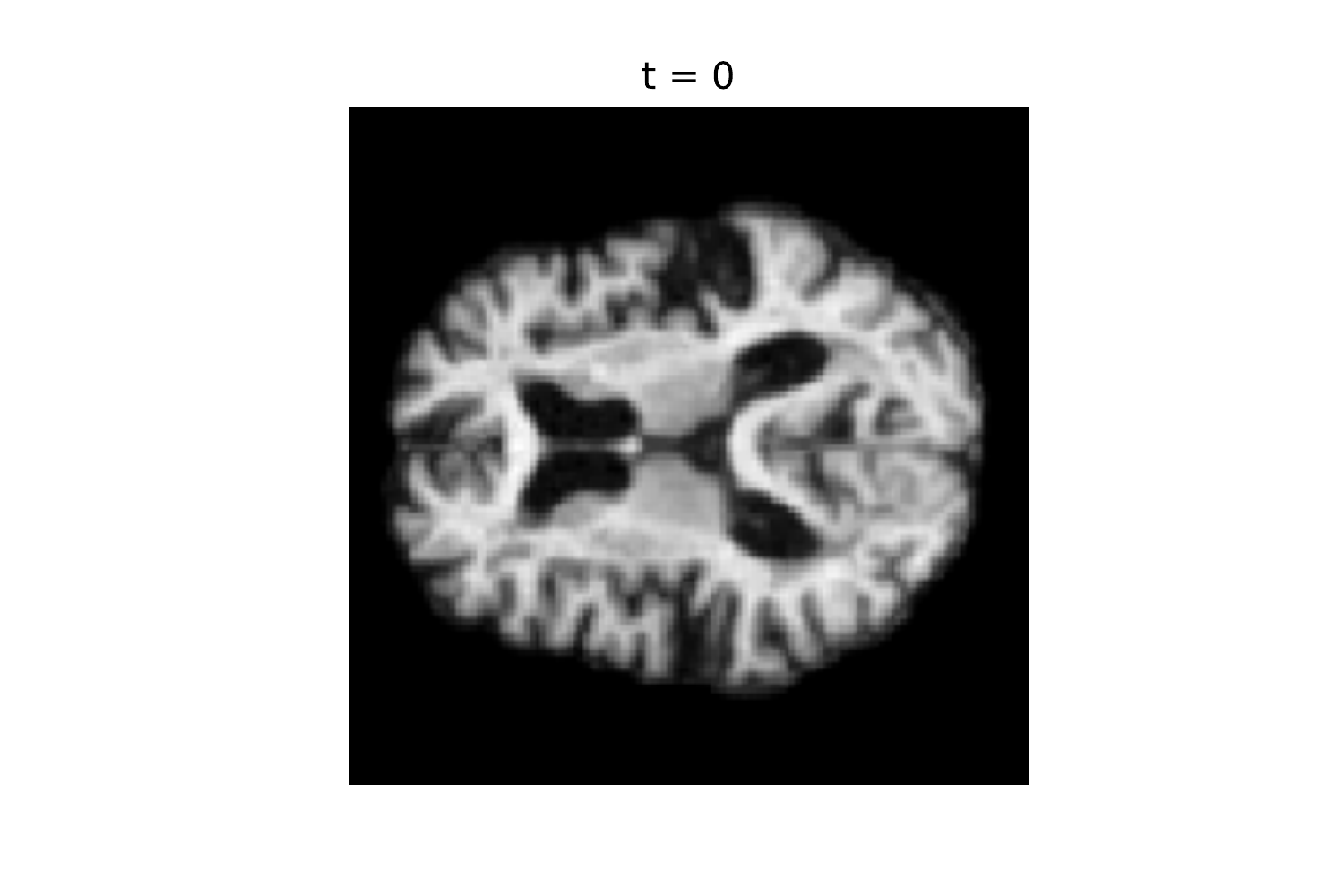}}
        \subfigure[$t=0.25$]{\includegraphics[width=0.24\columnwidth,trim=100 30 100 30,clip,angle=90]{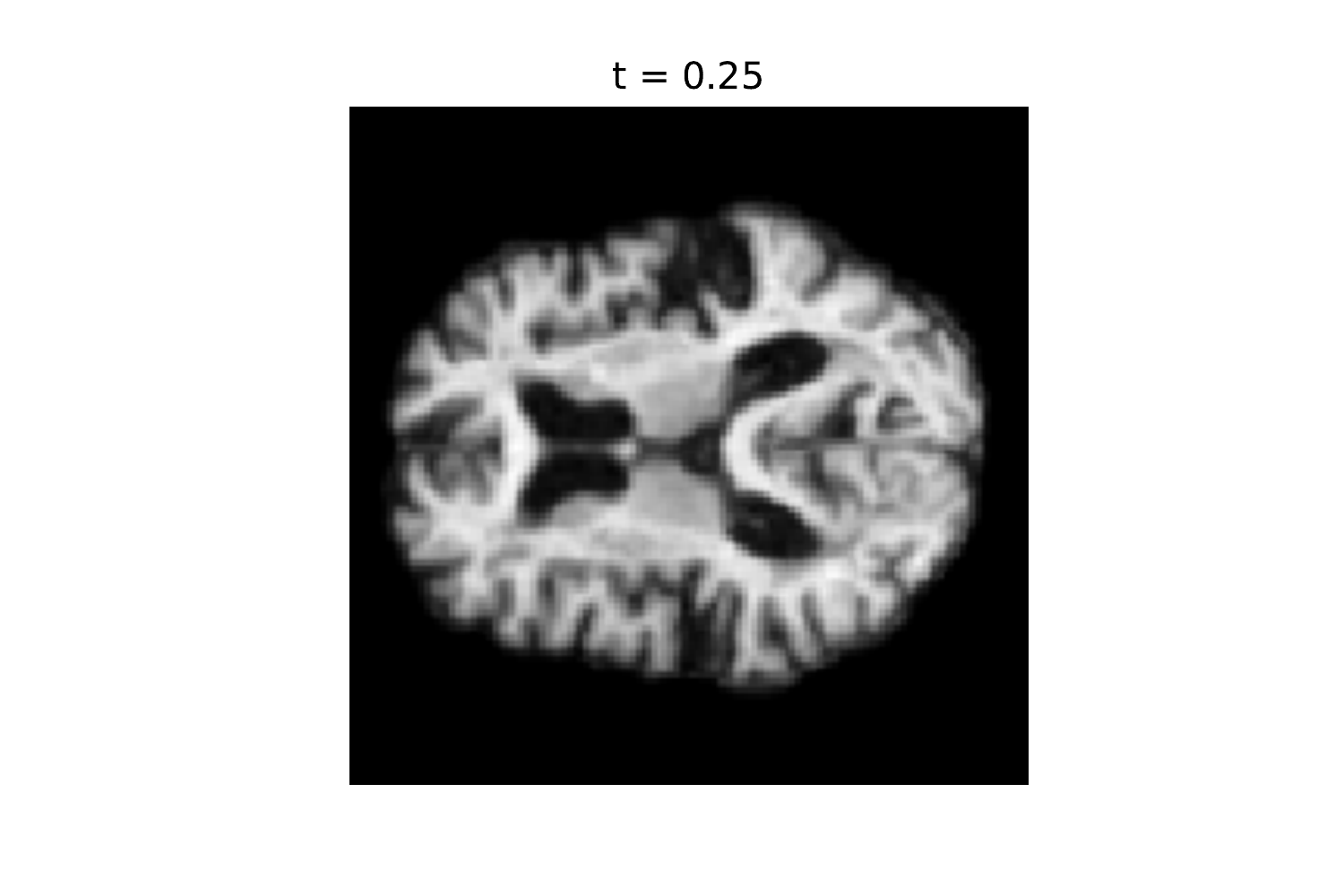}}
        \subfigure[$t=0.5$]{\includegraphics[width=0.24\columnwidth,trim=100 30 100 30,clip,angle=90]{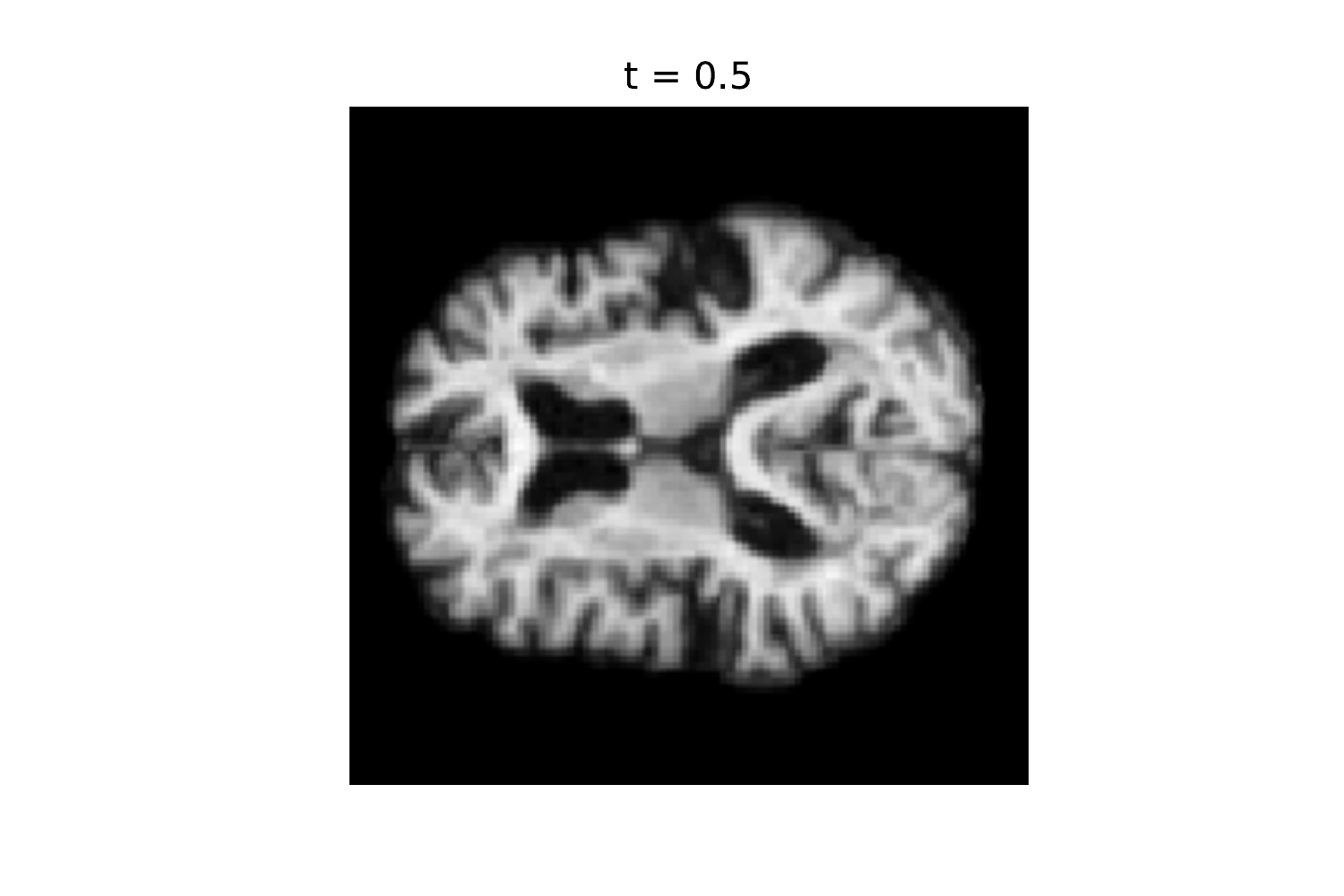}}
        \subfigure[$t=1$]{\includegraphics[width=0.24\columnwidth,trim=100 30 100 30,clip,angle=90]{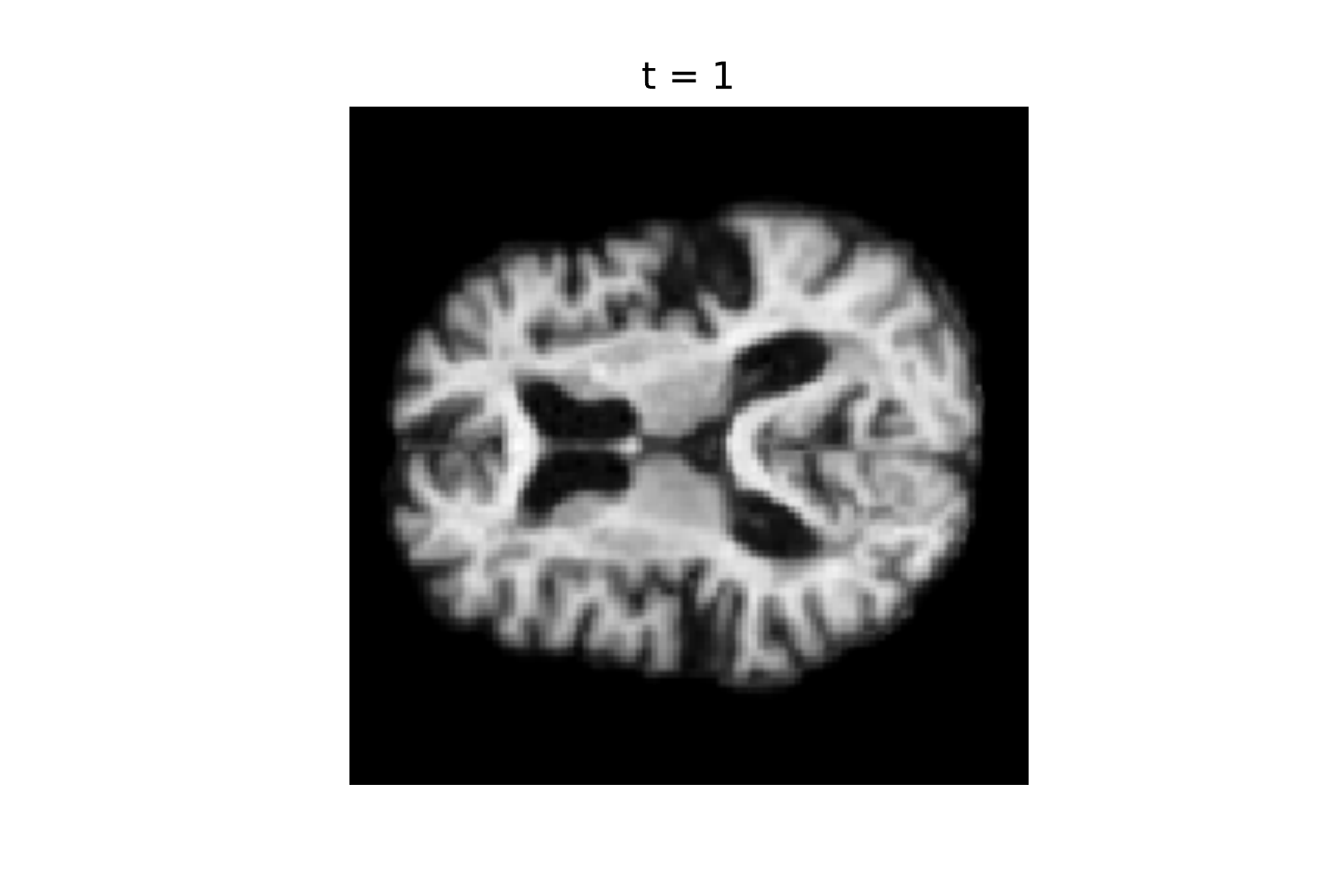}}
        \caption{4 time steps of a sample of the stochastic deformation process $I_0\circ\phi_t^{-1}$ at $t = 0,\, 0.25,\, 0.5,\, 1$.}
    \label{fig:stocev}
    \end{center}
\end{figure}

\vspace{-0.5cm}
\subsection{Inverse of the flow}

To calculate the moments of the stochastic deformation, we need to consider the nonlinear coupling between the process $\phi_t$ and the image. 
As the action of $\phi$ on $I_0$ is by composition with $\phi_t^{-1}$, we determine the SDE of the inverse $\phi_t^{-1}$, i.e. for each $t\in [0,1], \ \omega\in\Omega$, $\phi_t(\omega,\phi_t^{-1}(\omega,x,y)) = (x,y)$ for $(x,y)\in\mathcal{D}$. 
By the It\^o-Wentzel theorem~\cite{Krylov2011}, the SDE of $\phi_t^{-1}$ is,
\begin{align}
    d\phi_t^{-1} = -D\phi_t^{-1}v_t dt - \sum_{k=1}^p D\phi_t^{-1}\sigma_k\circ_S dW^k_t\, .
    \label{eq:phiinv}
\end{align}
By drawing sample paths of the stochastic flow $\phi_t^{-1}$, we obtain samples of deformed images under the model. 
Sample images at $t=1$ are shown in Fig. \ref{fig:init_var}, together with a plot of the generated image variation, and the noise fields used to simulate the sample data. 
Fig. \ref{fig:stocev} shows 4 time points from a sample path of the stochastic process $I_0\circ\phi_t^{-1}$ for $t=0,\, 0.25,\, 0.5,\, 1$. 

\subsection{Moments of Stochastic Image Deformation}\label{sec:taylor}

To approximate the first order moment of $\hat{I}_0\circ\phi_1^{-1}$, we consider a first order Taylor approximation of $\hat{I}_0\circ\phi_1^{-1}$ around the mean $\langle\phi_1^{-1}\rangle$, given as
\begin{align}
    \hat{I}_0\circ\phi_1^{-1}\approx \hat{I}_0\circ\langle\phi_1\rangle + \nabla(\hat{I}_0\circ\langle\phi_1^{-1}\rangle)^T(\phi_1^{-1} - \langle\phi_1^{-1}\rangle)\, .
\label{eq:taylor}
\end{align}
We consider two different means: $\langle\cdot\rangle$ which denote the mean of the stochastic processes $\phi_t$ and $v_t$, and $\mathbb{E}[\hat{I}_0\circ\phi_1^{-1}]$ for the expectation in image space.
Applying $\mathbb{E}$ to the Taylor approximation \eqref{eq:taylor} and using that $\nabla(\hat{I}_0\circ\langle\phi_1^{-1}\rangle)$ is deterministic and can be moved outside the mean function, we obtain the approximation
\begin{align}
    \mathbb{E}[\hat{I}_0\circ\phi_1^{-1}]\approx \hat{I}_0\circ\langle\phi_1^{-1}\rangle\, .
    \label{eq:mean}
\end{align}
In a similar manner, an approximation of the variance $\text{Var}[\hat{I}_0\circ\phi_1^{-1}]$ which singularly depend on the transition variance of the stochastic process $\phi_t^{-1}$, can be described by applying the first order Taylor approximation of $\hat{I}_0\circ\phi_1^{-1}$:
\begin{align}
    \text{Var}&\left[\hat{I}_0\circ\phi_1^{-1}\right] = \mathbb{E}\left[(\hat{I}_0\circ\phi_1^{-1})^2\right] - \left(\mathbb{E}[\hat{I}_0\circ\phi_1^{-1}]\right)^2 \nonumber \\
&\approx  \mathbb{E}\left[\left(\hat{I}_0\circ\langle\phi_1^{-1}\rangle + \left(\nabla(\hat{I}_0\circ\langle\phi_1^{-1}\rangle)\right)^T(\phi_1^{-1} -\langle\phi_1^{-1}\rangle)\right)^2\right] - \left(\mathbb{E}[\hat{I}_0\circ\phi_1^{-1}]\right)^2 \nonumber \\
&\approx  \left(\nabla(\hat{I}_0\circ\langle\phi_1^{-1}\rangle)^2\right)^T \left\langle\left(\phi_1^{-1} -\langle\phi_1^{-1}\rangle\right)^2\right\rangle\,.
\label{eq:var}
\end{align}
In the last approximation of \eqref{eq:var}, we used the approximation of the mean value presented in \eqref{eq:mean}.
To determine the moments of $\hat{I}_0\circ\phi_1^{-1}$, we therefore need the transition moments of the stochastic flow $\phi_t^{-1}$. These moments are studied in the next section by considering the Fokker-Planck equation.

\subsection{Moment Equations for $\phi_t^{-1}$ and $\boldsymbol{\tilde{v}}_t$}\label{sec:momphi}

 The transition distributions of the stochastic process $\phi_t^{-1}$ and the Fourier transformation of $v_t$ denoted $\boldsymbol{\tilde{v}}_t$, can be determined via the Kolmogorov forward equation also called the Fokker-Planck equation. 
 Based on this and the Kolmogorov operator $\mathcal{L}$, the moment evolution of a transformation of a stochastic process $X_t$ under a real-valued function $h$ is described by
\begin{align}
    \frac{d}{dt}\langle h(X_t)\rangle = \langle\mathcal{L}h(X_t)\rangle \, .
\label{eq:mom}
\end{align}
 More information on the procedure for calculating the moment evolution via the Fokker-Planck equation can be found in~\cite{arnaudon_geometric_2017}. 

For sake of notation, let $\psi_t = \phi_t^{-1}$ and $\psi_{ij} = \psi_t(x_i,y_j) = (\psi_t^x(x_i,y_j),\psi_t^y(x_i,y_j))$ for a discretisation of the image domain $\mathcal{D}$ to a grid $(x_i,y_j)_{ij}$. 
When considering a stochastic process, e.g. $\psi_t$, described by an It\^o SDE, the Kolmogorov operator $\mathcal{L}$ is found by applying It\^o's formula. 
The Kolmogorov operator is defined as the drift of the resulting stochastic process. 
As an example, the Kolmogorov operator for the 2-dimensional stochastic process $\psi_{ij,t}$ is given by 
\begin{align*}
    \mathcal{L}f &= \frac{\partial f}{\partial t} + \sum_{\alpha\in\{x,y\}} \Bigg[\frac{\partial f}{\partial \psi_{ij}^\alpha}(-D\psi_{ij}v_t(x_i,y_j) \\
& + \frac{1}{2}\sum_{k}D[D\psi_{ij}\sigma_k(x_i,y_j)]D\psi_{ij}\sigma_k(x_i,y_j))_\alpha\Bigg]
+ \frac{1}{2}\sum_{\alpha,\beta\in\{x,y\}} C_{\alpha\beta}\frac{\partial^2 f}{\partial \psi_{ij}^\alpha\partial \psi_{ij}^\beta} \, , 
\end{align*}
for $C = bb^T$, $b = (D\psi_{ij}\sigma_1(x_i,y_j)\cdots D\psi_{ij}\sigma_p(x_i,y_j))$, and a twice differentiable real-valued function $f$.

For the Kolmogorov operator for $\psi_{ij}$, the evolution of the moments are determined by \eqref{eq:mom}. Let $\gamma\in\{x,y\}$ be given and $(x_i,y_j)$ be a fixed pixel in the grid. 
Consider the time evolution of the first moment of $\psi_{ij}^\gamma$
\begin{align}
    \frac{\partial}{\partial t}&\langle\psi_{ij}^\gamma\rangle = \langle\mathcal{L}\psi_{ij}^\gamma\rangle  \nonumber \\
    &=\left\langle(-D\psi_{ij}v_t(x_i,y_j) + \frac{1}{2}\sum_{k}D[D\psi_{ij}\sigma_k(x_i,y_j)]D\psi_{ij}\sigma_k(x_i,y_j))_\gamma\right\rangle \label{eq:approx} \\
    &\approx (-D\langle\psi_{ij}\rangle\langle v_t(x_i,y_j)\rangle + \frac{1}{2}\sum_{k}D[D\langle\psi_{ij}\rangle\sigma_k(x_i,y_j)]D\langle\psi_{ij}\rangle\sigma_k(x_i,y_j))_\gamma\, . \nonumber
\end{align}
Notice the coarse approximation used in \eqref{eq:approx} to describe the time evolution uniquely by moments of $\psi_t$ and $v_t$. 
This approximation is used in the rest of the paper and assumes independence of the random variables. 
That is, any mean of a product of random variables $\psi^\gamma_{ij}$, $v^\gamma_{ij}$ is approximated by the product of the first order moments of each random variable. 
Extending the equations to include higher order correlation will be the topic of future works. 

The derivation of the variance of $\psi_{ij}^\gamma$ is calculated as above and results in the moment evolution
\begin{align*}
    &\frac{\partial}{\partial t}\left\langle(\psi_{ij}^\gamma - \langle\psi_{ij}^\gamma\right\rangle)^2\rangle = \langle\mathcal{L}(\psi_{ij}^\gamma - \langle\psi_{ij}^\gamma\rangle)^2\rangle \\
&= \Bigg\langle 2(\psi_{ij}^\gamma - \langle\psi_{ij}^\gamma\rangle) (-D\psi_{ij}v_t(x_i,y_j) +\frac{1}{2}\sum_{k}D(D\psi_{ij}\sigma_k(x_i,y_j))D\psi_{ij}\sigma_k(x_i,y_j))^\gamma\Bigg\rangle \\
&\hspace{1cm}+\langle C_{\gamma\gamma}\rangle \approx \langle C_{\gamma\gamma}\rangle\, ,
\end{align*}
where $C$ is given as above and where the assumption of independence between variables is used to split the first term into the product of $\langle\psi_{ij}^\gamma - \langle\psi_{ij}^\gamma\rangle\rangle$ and $\left\langle -D\psi_{ij}v_t(x_i,y_j)+\frac{1}{2}\sum_{k}D(D\psi_{ij}\sigma_k(x_i,y_j))D\psi_{ij}\sigma_k(x_i,y_j))^\gamma\right\rangle$.

The moment evolution of $\psi_{ij}^\gamma$ in \eqref{eq:approx} depends on the first order moment of the time-varying velocity field $v_t$. 
As described in Section \ref{sec:stocev}, applying a spatial discrete Fourier transform to a discretization, $\boldsymbol{v}_t$, of $v_t$ to a grid $(x_i,y_j)_{ij}$, results in a computationally feasible optimisation procedure. 
Using the property $\mathcal{F}[\langle \boldsymbol{v}_t^\alpha\rangle] = \langle\mathcal{F}[\boldsymbol{v}_t^\alpha]\rangle = \langle\boldsymbol{\tilde{v}}_t\rangle$, the moment evolution of $\boldsymbol{\tilde{v}}_t$ is calculated by the same procedure as above applied to the Fourier transform of the It\^o\ SDE presented in Section~\ref{sec:mm}. 
The moment evolution is given as
\begin{align*}
    \frac{d}{dt}\langle\tilde{v}_{t,ij}^\alpha\rangle &= \langle\mathcal{L}\tilde{v}_{t,ij}^\alpha\rangle  \\
    \approx &-\left(K\text{ad}^*_{\langle\tilde{v}_t\rangle}\langle\tilde{m}_t\rangle\right)_{ij}^\alpha + \frac{1}{2}\sum_{k=1}^p\left((\tilde{D}[K\text{ad}^*_{\tilde{\sigma}_k}\langle\tilde{m}_t\rangle)*(K\text{ad}^*_{\tilde{\sigma}_k}\langle\tilde{m}_t\rangle)\right)_{ij}^\alpha.
\end{align*}
The above moment equation is described in Fourier domain. Here, $*$ denotes convolution, and $\tilde{D}$ a central difference Jacobian matrix (for more information see \cite{zhang2015finite}).
Using the Fourier representaition, it takes approximately $8$ seconds to solve the moment equations in a situation of 1 noise field, 100 time steps, truncated to 16 Fourier frequencies and on a standard laptop.

\subsection{Similarity of Moment Images}

With the method of moments, we seek to maximise the similarity between the moments of the observed data and moments of the stochastic deformation of the initial image $\hat{I}_0$. 
For an example of the variance of the deformed images and the data sample variance, see Fig. \ref{fig:res}. Note that the approximation of the mean $\hat{I}_0\circ\langle\phi_1^{-1}\rangle$ and the variance $\text{Var}[\hat{I}_0\circ\phi_1^{-1}]$ are images themselves. Therefore, moment matching turns into matching of images, however not a match of observed images as regularly performed in image registration, but instead a match of $\hat{I}_0\circ\langle\phi_1^1\rangle$ and $\text{Var}[\hat{I}_0\circ\phi_1^{-1}]$ towards their sample equivalents. Interestingly, we see that classical image similarity measures can be used to compare these images effectively.

Due to the approximations described in Section \ref{sec:taylor}, we cannot expect a perfect match between moment images of model and data samples. Because of this, we find that the $L^2$-distance often does not result in good matching. Instead, using normalised mutual information, we were able to retrieve the correct values of the variance parameters. 
However, as normalised mutual information corrects for intensity differences, it is generally invariant to changes in the noise amplitude parameters $\lambda_k$. 
Therefore, we use normalised mutual information to get a good match for the variance parameters before estimating the amplitude parameters using a combination of $L^2$ distance and unnormalised mutual information.

We let the noise kernels be Gaussian 
$
\sigma_k(\boldsymbol{x}) = \lambda_k\text{exp}\left(\frac{\|\boldsymbol{x}-\mu_k\|^2}{2\tau^2_k}\right)\mathrm{Id}_2
$ for simplicity.
Hence, the parameters to be estimated consist of the amplitudes $\lambda_k\in\mathbb{R}$ and the variances $\tau^2_k$. 
Varying the variances $\tau^2_k$ changes the spatial effect of the noise fields while varying the amplitude $\lambda_k$ affects the intensity of the noise.
Let $\mu_1$ denote the sample average, i.e. $\mu_1 = \frac{1}{n}\sum_{i=1}^n I^i_1$. 
The optimisation procedure first optimise the function,
\begin{align}
    f(\{\tau_k\}_{k=1,\ldots,p}) =\, & \text{MI}_{\text{norm}}\left(\hat{I}_0\circ\langle\phi_1^1\rangle_{(\tau_k,\lambda_k)},\ \mu_1\right) \nonumber \\ &+ \text{MI}_{\text{norm}}\left(\text{Var}_{(\tau_k,\lambda_k)}[\hat{I}_0\circ\phi_1^{-1}],\ \frac{1}{n}\sum_{i=1}^n(I^i_1-\mu_1)^2\right)\, , 
\end{align}
for the variance parameters $\tau^2_k$ and then optimise the objective function 
\begin{align}
    g(\{\lambda_k\}_{k=1,\ldots,p}) =  &\left\|\hat{I}_0\circ\langle\phi_1^1\rangle_{(\tau_k,\lambda_k)}-\mu_1\right\| \nonumber + \text{MI}\left(\hat{I}_0\circ\langle\phi_1^1\rangle_{(\tau_k,\lambda_k)},\ \mu_1\right) \\
&+ \left\|\text{Var}_{(\tau_k,\lambda_k)}[\hat{I}_0\circ\phi_1^{-1}]-\frac{1}{n}\sum_{i=1}^n(I^i_1-\mu_1)^2\right\| \nonumber \\
&+ \text{MI}\left(\text{Var}_{(\tau_k,\lambda_k)}[\hat{I}_0\circ\phi_1^{-1}],\ \frac{1}{n}\sum_{i=1}^n(I^i_1-\mu_1)^2\right) \, , 
\label{eq:objlamb}
\end{align}
for the amplitude $\lambda_k$. 

\section{Simulation Study}\label{sec:exp}

In this section, we present a simulation study aiming at illustrating the ability of the framework to infer parameters given MR brain images.
Given an initial 2D MRI image from the OASIS database, two datasets of 200 observations were simulated based on $9$ noise fields located in a grid on the image domain.

\begin{figure}[htpb]
    \begin{center}
        \subfigure{\includegraphics[width=0.245\columnwidth,trim=100 30 100 30,clip,angle=90]{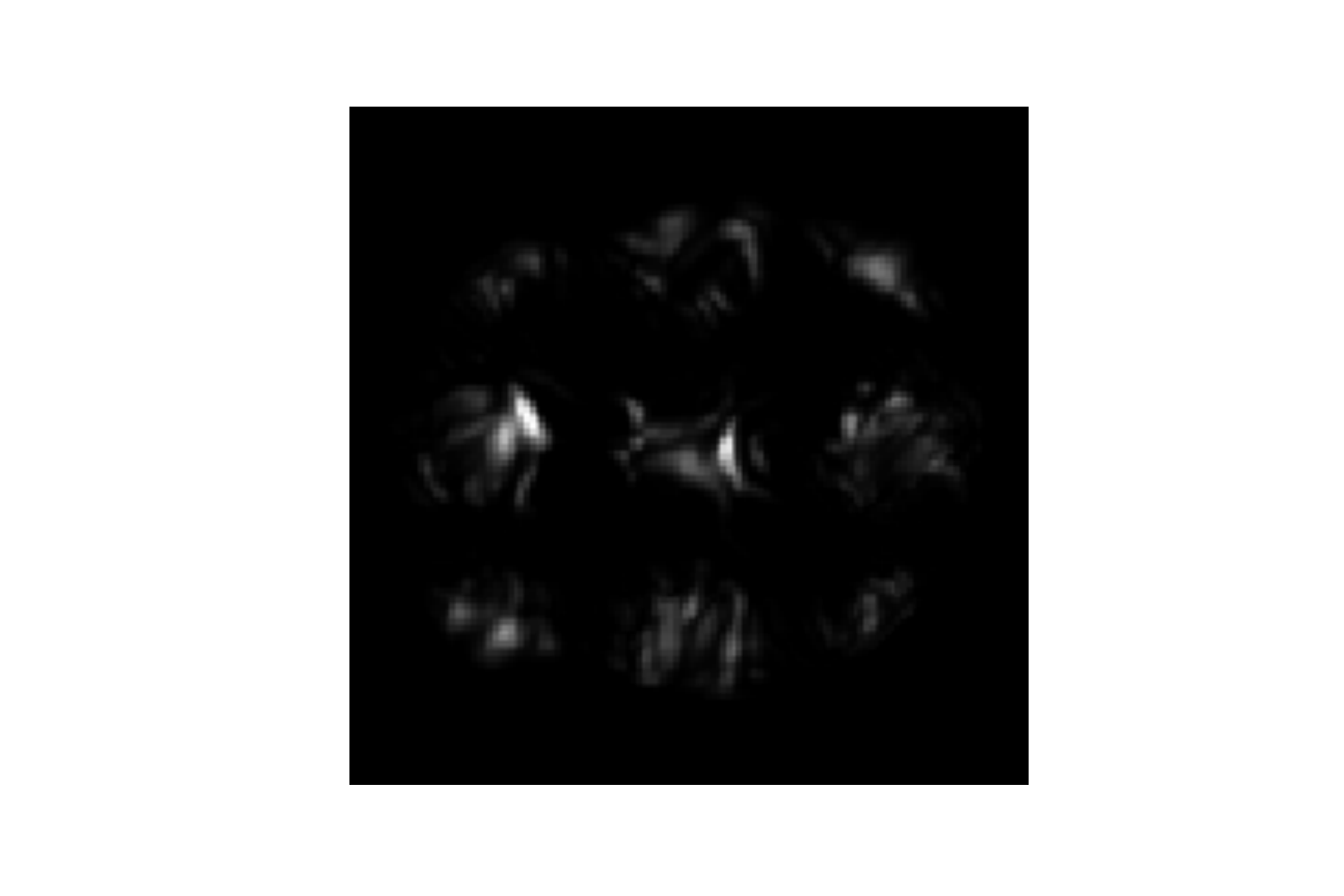}}
        \subfigure{\includegraphics[width=0.245\columnwidth,trim=100 30 100 30,clip,angle=90]{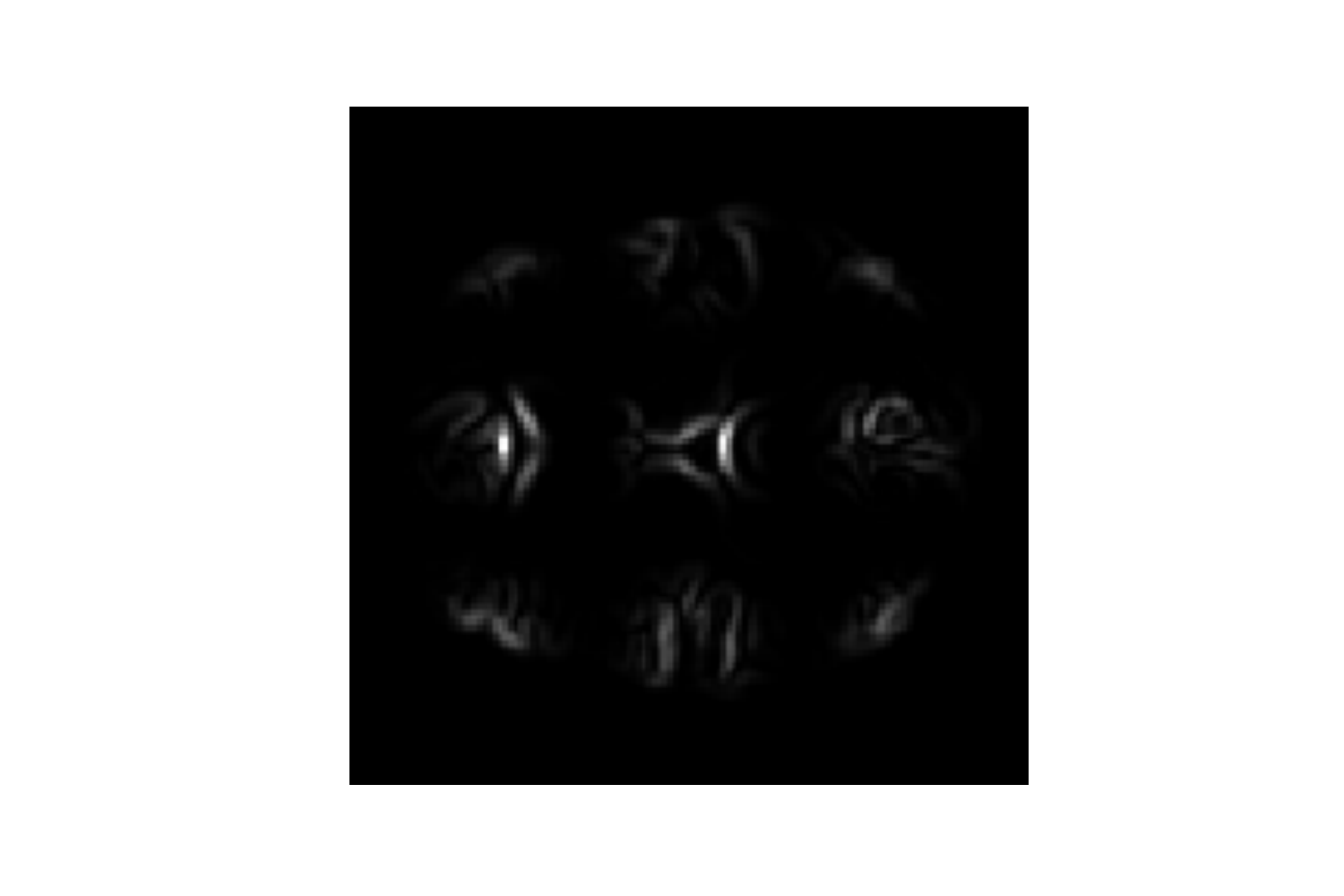}}%
        \subfigure{\includegraphics[width=0.245\columnwidth,trim=100 30 100 30,clip,angle=90]{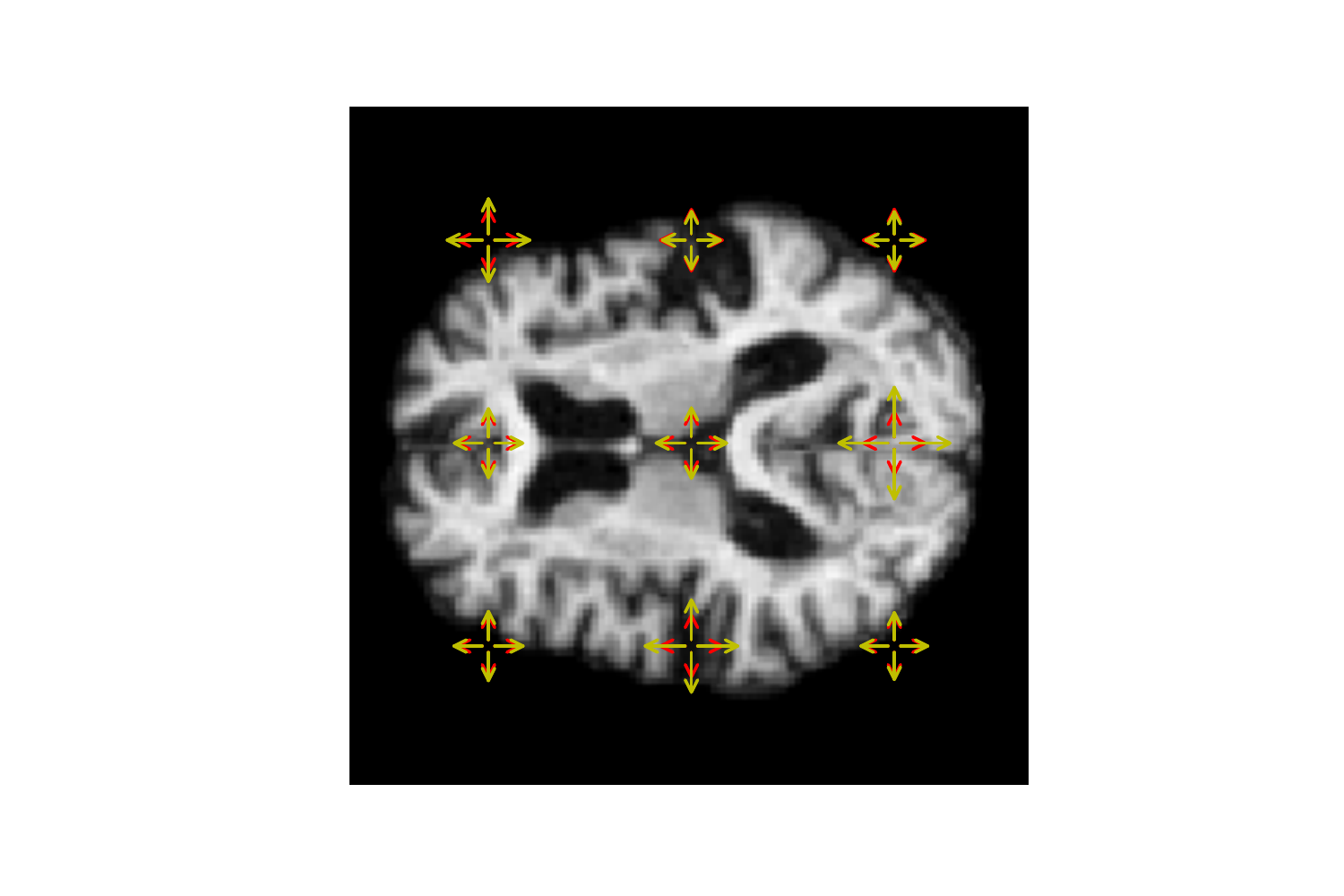}}%
        \subfigure{\includegraphics[width=0.245\columnwidth,trim=100 30 100 30,clip,angle=90]{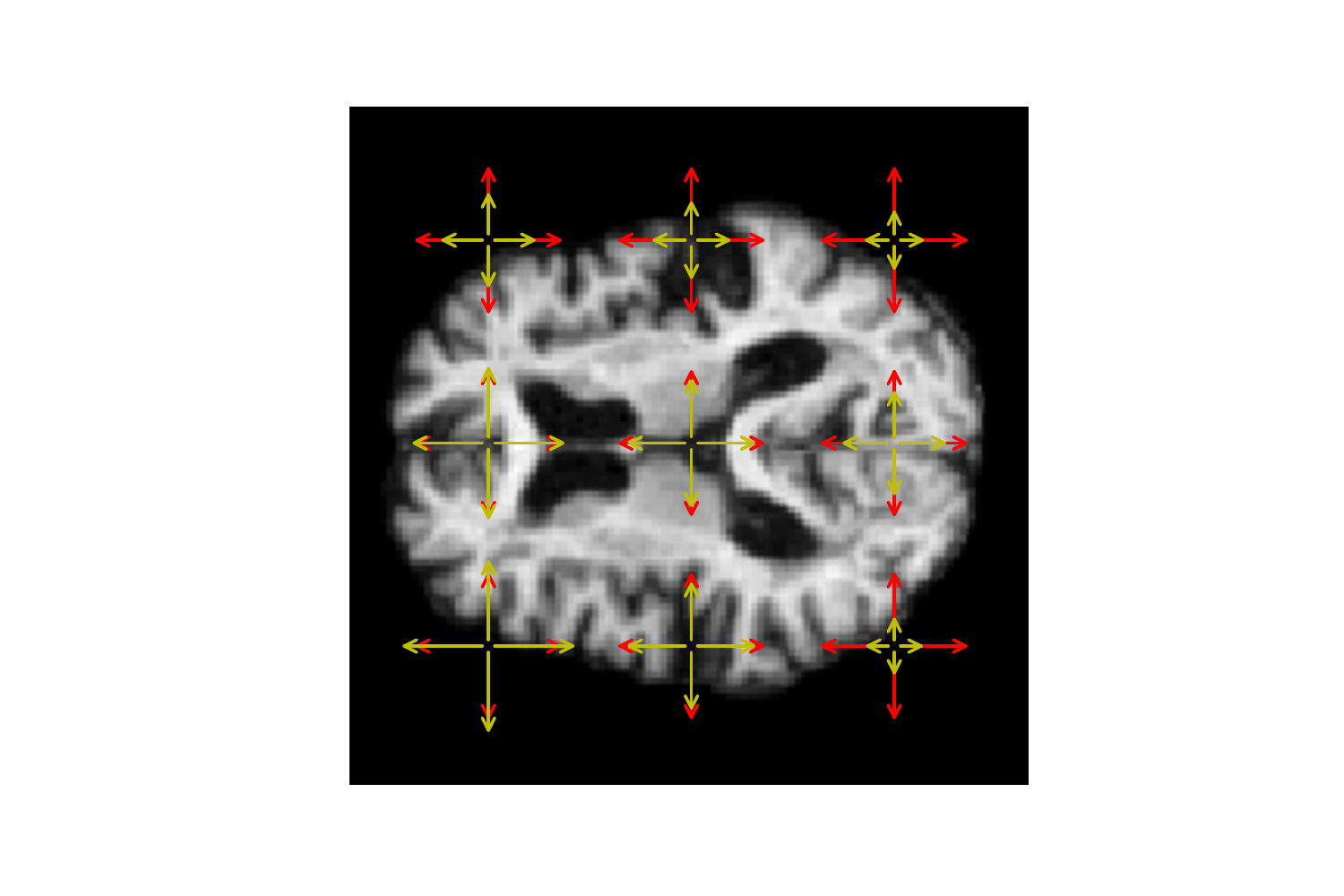}}
        \caption{\textit{From left:} (1.) Pixel-wise variation in the data. (2.) Variation estimated by the model. (3.) Estimates of the variance parameter $\tau_k^2$ compaired to the true values (4.) Estimates of the variance and amplitude parameter compaired to the true values. The length of the arrows correspond to $\tau_k\lambda_k$. The arrows show the location and width of the noise fields. The red arrows correspond to the true values, while the yellow defines the resulting estimated parameters.} 
    \label{fig:res}
    \end{center}
\end{figure}

The first data sample was simulated based on 9 Gaussian noise fields. The location of the noise fields is shown in Fig. \ref{fig:init_var}, which also shows the initial brain image, $I_0$. 
The variation of the simulated data sample is visualized in Fig. \ref{fig:res}.

In the simulation study, estimates of the variance parameters $\tau_k^2$ and amplitude $\lambda_k$ were found for the $9$ noise fields. 
The true value of the standard deviation $\tau_k$ was set to $0.06$ and the amplitudes $\lambda_k = 2.0$ for all $k=1,\ldots,9$. 

The initial values of the optimisation procedure were found by investigating the parameter space randomly $40$ times picking the values with the smallest objective value. 
A gradient descent estimation was used with a line search for determining the step size at each iteration.

Fig.~\ref{fig:res} shows the initial brain, $I_0$, with a comparison of the true values of $\tau_k$ and $\lambda_k$, and the values found by the optimisation procedure. The red arrows are the true values, and the yellow defines the estimated parameter values.
 The model is able to retrieve the parameters of $\tau_k$ for all $k=1,\ldots,9$. It also returns a good estimate of the amplitude parameters, in particular for the noise fields located inside the brain. For noise fields on the boundary of the brain or in the background, the model does not have access to enough information in the intensity differences to return precise estimates of the amplitude parameters.  

\begin{figure}[htpb]
    \begin{center}
        \subfigure{\includegraphics[width=0.33\columnwidth,trim=80 60 40 70,clip]{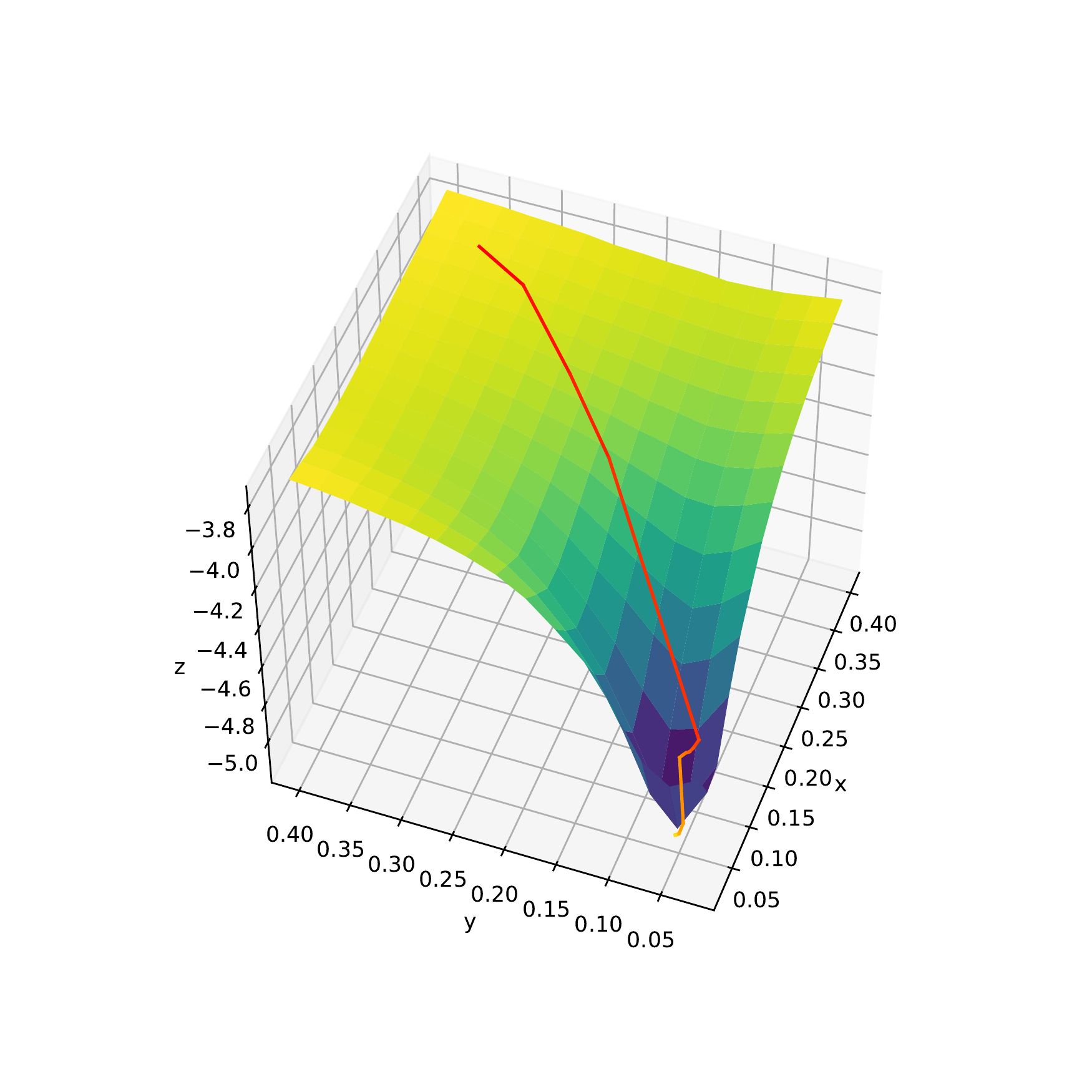}}%
        \subfigure{\includegraphics[width=0.33\columnwidth,trim=80 60 40 70,clip]{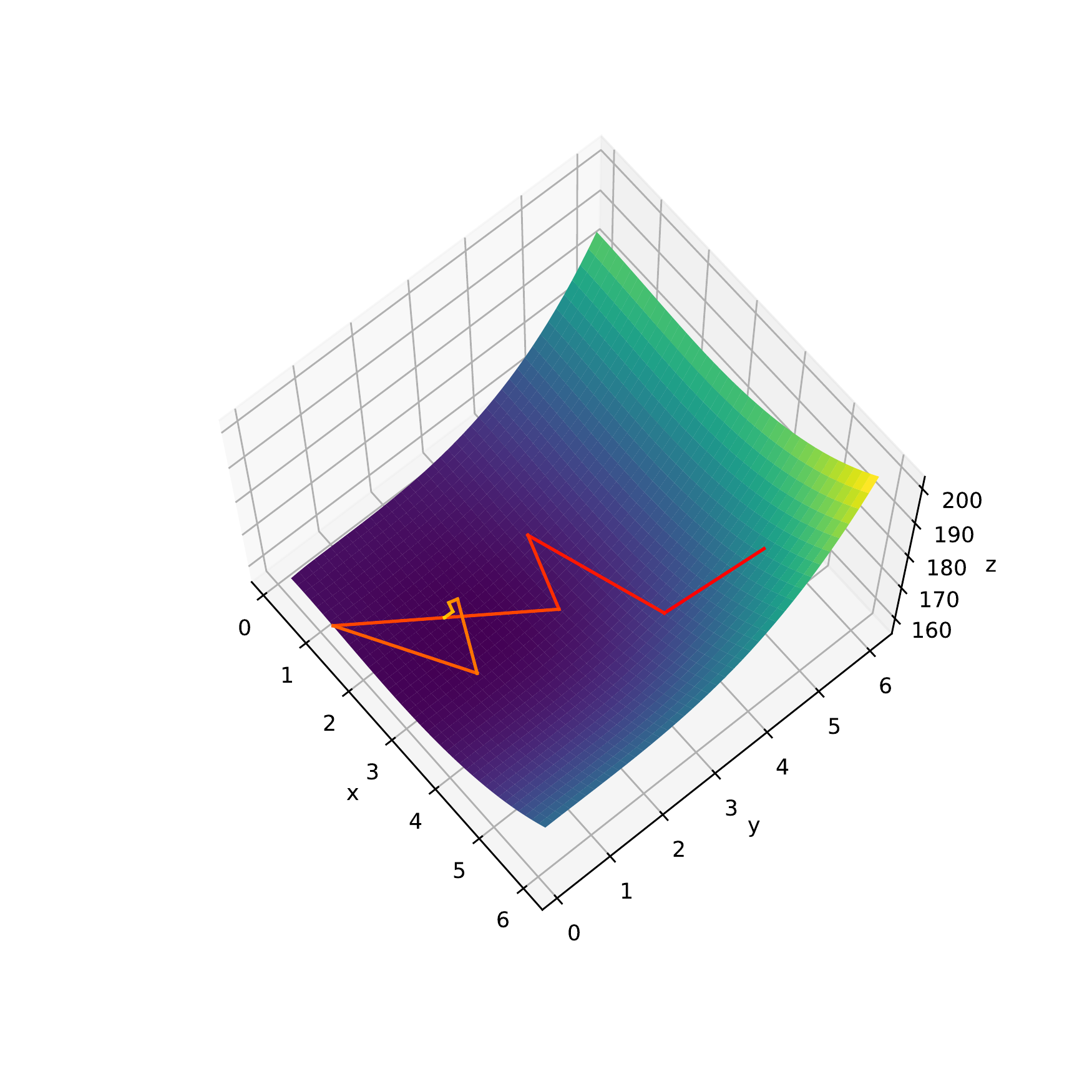}}%
        \subfigure{\includegraphics[width=0.33\columnwidth,trim=80 60 40 70,clip]{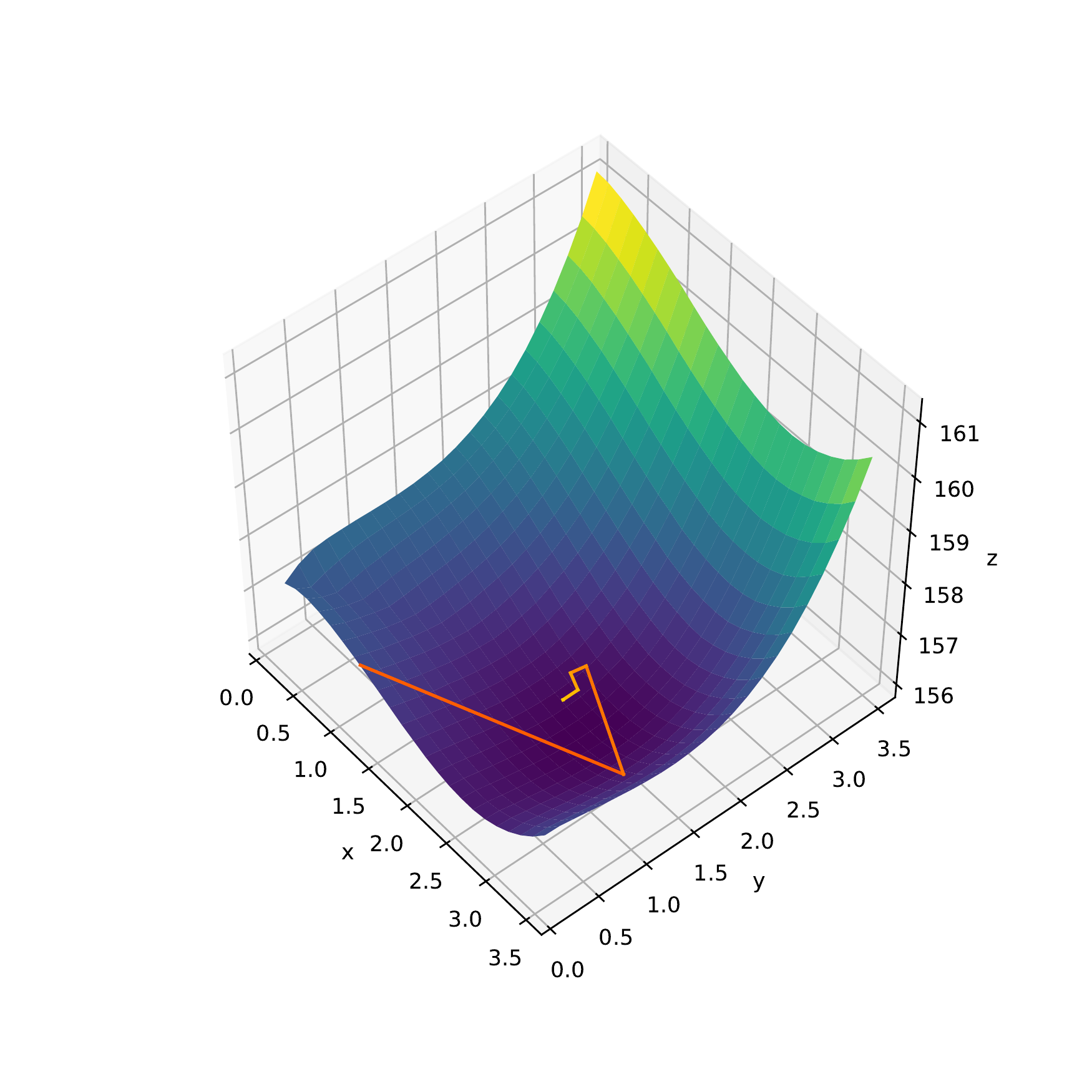}}
        \caption{\textit{From left:} (1.) Objective function for the optimisation of variance parameters $\tau_k^2$ for an example of $2$ noise fields with fixed amplitude parameters. (2.) The objective function for the amplitude $\lambda_k$ for an example of $2$ noise fields with fixed variance parameters. (3.) A zoom of the previous image (2.) around the estimated value.}
    \label{fig:objective}
    \end{center}
\end{figure}

To give an intuition of convergence of the optimisation, Fig.~\ref{fig:objective} (left) shows the objective function with gradient steps for the optimisation of $\tau_k$ in the case of $k=2$. The amplitude is held fixed in this situation. In the same figure is shown the objective function of the amplitude $\lambda_k$ when the variance is held fixed.

\begin{figure}[htpb]
    \begin{center}
        \subfigure{\includegraphics[width=0.245\columnwidth,trim=100 33 100 30,clip,angle=90]{Figures/9noise/data_var901.pdf}}%
        \subfigure{\includegraphics[width=0.245\columnwidth,trim=100 33 100 30,clip,angle=90]{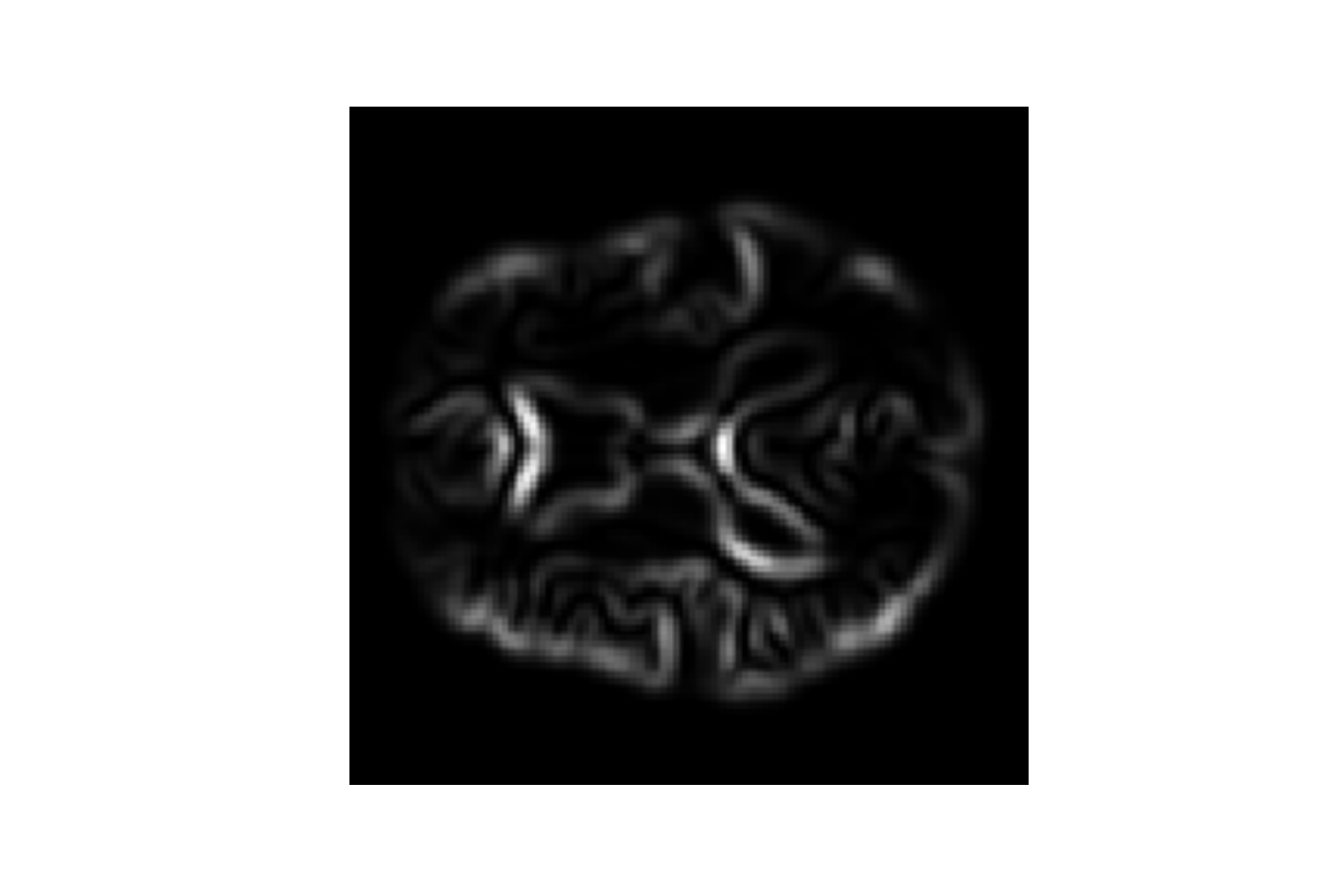}}
        \subfigure{\includegraphics[width=0.245\columnwidth,trim=100 33 100 30,clip,angle=90]{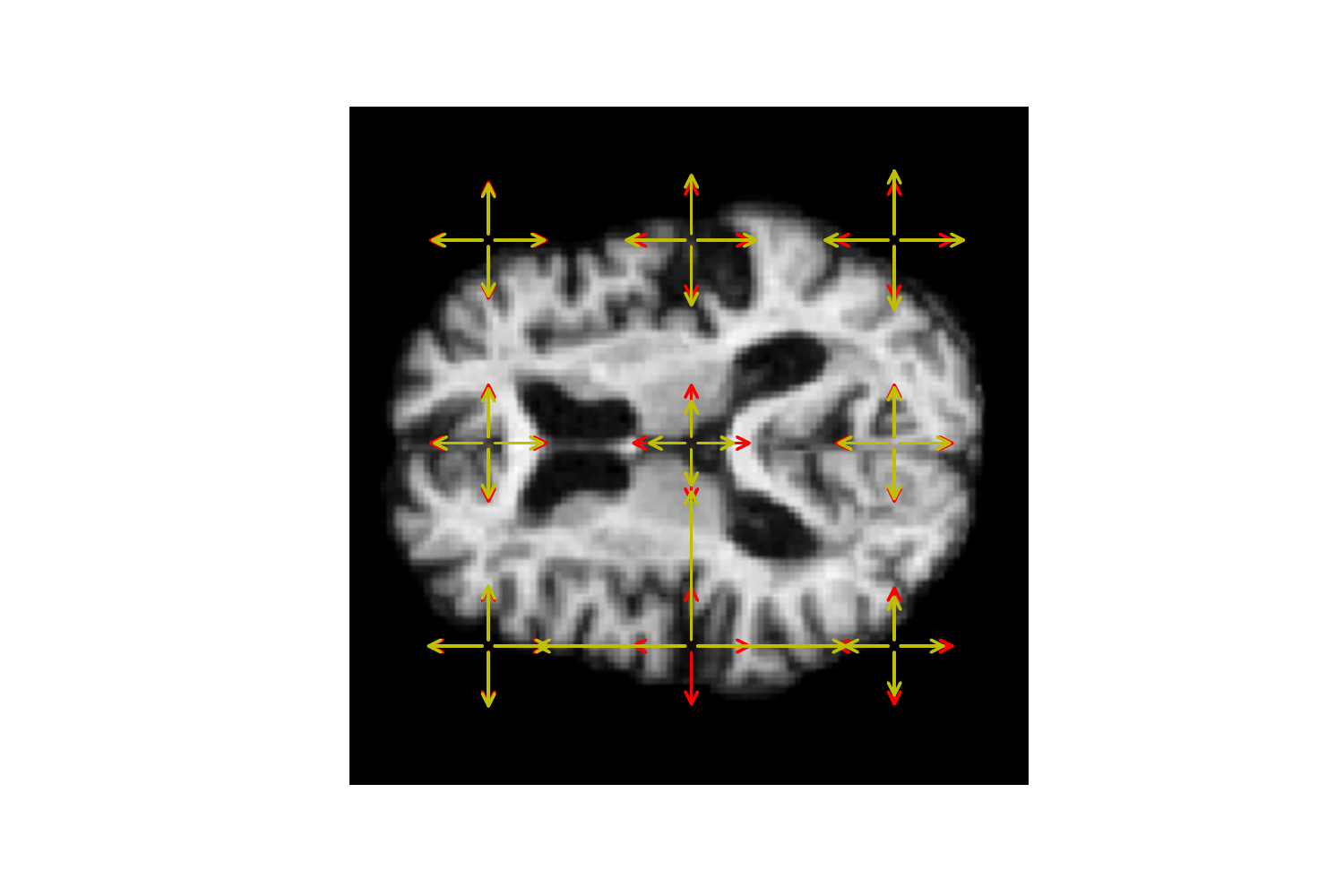}}
        \caption{(left) Data variation. (middle) estimated variation from the model. (right) Estimated variance parameters. The red arrows correspond to the true values, the yellow defines the resulting estimated parameters.}
    \label{fig:res01}
    \end{center}
\end{figure}

The second dataset was simulated based on $9$ noise fields with larger deviation at $\tau_k = 0.1$ and where the area of variation of each field intersect. The $9$ noise fields are shown in Fig. \ref{fig:init_var}. The amplitude of the noise fields is again set to $\lambda_k = 2.0$. Retrieving the true values of the parameter vector is harder in this case as the optimisation procedure is more prone to reach a local minimum. Therefore, we focus this example on estimating the variance parameters $\tau_k^2$. As shown in Fig. \ref{fig:res01}, good estimates for the variance parameters of most of the noise fields are obtained.
Estimation of both variance and amplitude with large intersection between the fields is challenging because little information in the data is availble to precisely determine to which of the intersecting noise fields observed variation belong. This could be handled either by imposing a prior on the noise to enforce spatial regularity of the noise amplitudes or by considering noise that is naturally uniform over the domain, e.g. by representing the noise itself in Fourier domain.

\vspace{-0.5cm}
\section{Conclusion}
\label{sec:con}

We presented a model for estimating the variation in medical images occurring from time-continuous deformation variation. 
The model was based on the stochastic generalisation of the LDDMM framework, using the FLASH procedure to make a natural dimensionality reduction resulting in computationally fast image deformations. 
Determining the moments of the stochastic flow of the deformations $\phi_t$ and velocity fields $v_t$, the method of moments was applied to estimate the parameter vector for noise fields defining the variation in the data sample. 
The moments of $\phi_t$ and $v_t$ have been calculated using the Fokker-Planck equation of the evolution of the truncated Fourier expansion of the image deformation. 
These moments were compared to the data distribution allowing for parameter estimation.

For future work, a natural extension is to define the noise fields in the Fourier domain instead of the spatial fields discussed here.

To calculate the image moments, a coarse Taylor approximation was applied. In future work, we wish to perform a broader investigation of the consequence of making this approximation and whether alternative methods can be used to get a more precise estimation of the image moments.

We disregarded the subject-specific variation and initialised the stochastic deformation by a single image. However, in real data, the subject-specific variation can be large and the model will generally not return a good estimate of data variation when this effect is not taken into account. The model can be extended to include a subject-specific initial image, such that only the variation over time for each individual is modelled and not the total population variation.

Finally, the model has been applied to 2D slices of 3D images. The model is fully general, and since the FLASH framework can handle 3D image data, we wish to include analyses of 3D images in the stochastic framework.
\\

\textbf{Acknowledgements.} AA acknowledges EPSRC funding through award EP/N014529/1 via
the EPSRC Centre for Mathematics of Precision Healthcare. LK and SS are supported by the CSGB
Centre for Stochastic Geometry and Advanced Bioimaging funded by a grant from the
Villum Foundation. 

\bibliographystyle{plain}

\end{document}